\documentclass[a4paper]{amsart}

\usepackage{amsmath,latexsym}
\usepackage{amssymb}
\usepackage{graphics, graphicx}
\usepackage{amsfonts}
\usepackage[nesting]{hyperref}
\usepackage{pdflscape}
\usepackage[justification=centering]{caption}
\usepackage{subcaption}

\newcommand{\dsum}{\displaystyle\sum}

\newcommand{\dmax}{\displaystyle\max}

\def\R{\mathbb{R}}
\def\N{\mathbb{N}}

\def\k{\mathrm{k}}
\def\sign{\mathrm{sign}}

\def\r{\texttt{R}}
\def\flp{f_{\lambda,p}}
\def\R{\mathbb{R}}
\def\Z{\mathbb{Z}}

\def\ext{\mathrm{Ext}}
\def\Q{\mathbb{Q}}
\def\N{\mathbb{N}}

\def\bbeta{\boldsymbol\beta}
\def\bepsilon{\boldsymbol\varepsilon}

\def\GCoD{\mathrm{GCoD}}
\def\d{\mathrm{D}}
\def\H{\mathcal{H}}

\usepackage{multicol}
\usepackage{multirow}
\usepackage[english, activeacute]{babel}
\usepackage[latin1]{inputenc}

\usepackage{hyperref}
\usepackage{color}

\usepackage[margin=2cm]{geometry}

\newtheorem{theo}{Theorem}

\newtheorem{ex}[theo]{Example}
\newtheorem{lem}[theo]{Lemma}
\newtheorem{prop}[theo]{Proposition}
\newtheorem{cor}[theo]{Corollary}
\newtheorem{rem}[theo]{Remark}

\usepackage{multicol}

 \definecolor{dblue}{rgb}{0.2, 0.2, 0.6}
 \definecolor{dred}{rgb}{0.75, 0.0, 0.2}

\usepackage{pgfplots}
\pgfplotsset{compat=newest}
\usetikzlibrary{decorations.markings}    % Arrows

\begin{document}

\title{A general framework for locating hyperplanes to fitting set of points}

\author{V\'ictor Blanco}

\address{Dept. Quantitative Methods for Economics \& Business, Universidad de Granada.}
\email{vblanco@ugr.es}

\author{Justo Puerto}
\address{Dept. Estad{\'\i}stica e Investigaci\'on Operativa, Universidad de Sevilla.}
\email{puerto@us.es}

\author{Rom\'an Salmer\'on}
\address{Dept. Quantitative Methods for Economics \& Business, Universidad de Granada.}
\email{romansg@ugr.es}



\begin{abstract}
This paper presents a  family of new methods for locating/fitting hyperplanes with respect to a given set of points. We introduce a general framework for a family of  aggregation criteria of different  distance-based errors. The most popular methods found in the specialized literature can be cast within  this family as particular choices of the errors and the aggregation criteria. Mathematical programming formulations for these methods are stated and some interesting cases are analyzed. It is also proposed a   new goodness of fitting index which extends the classical coefficient of determination. A series of illustrative examples and extensive computational experiments implemented in \r{} are provided to show the performances of some of the proposed methods.
\end{abstract}

\maketitle

\section{Introduction}

The problem of locating hyperplanes with respect to a given set of point is well-known in Location Analysis \cite{schobelbook}. This problem is closely related to another common question in Data Analysis: to study the behavior of a given set of data with respect to a fitting body expressed with an equation of the form $f(X_1, \ldots, X_d)=0$. This last problem reduces to the estimation of the `best' function $f$ that expresses the relationship between the provided data or in other words to the location of the surface $f(X)=0$ that minimizes some aggregation function of the distances of the points (data set) to the dimensional facility $f$ (see \cite{ACT16,DMS,DSW02}). In many cases and for the sake of simplicity, the family of functions where $f$ belongs to is usually fixed and then, real  parameters of such a function must be determined. The most widely used family of functions considered in this framework, probably because of its simplicity, is the family of linear functions, namely the above equation is of the form $f(X_1, \ldots, X_d)=\beta_0 + \sum_{k=1}^d \beta_k \, X_k = 0$ for $\beta_0, \beta_1, \ldots, \beta_d \in \R$.

To perform such a fitting, we are given a set of points $\{x_1, \ldots, x_n\} \subset \R^d$, and one tries to find the values $\hat \bbeta = (\hat \beta_0, \hat \beta_1, \ldots, \hat \beta_d)$ that minimize some measure of the deviation  of the data with respect to the hyperplane  $\mathcal{H}(\hat \bbeta) = \{z \in \R^d: \hat \beta_0 + \sum_{k=1}^d \hat \beta_k z_k=0\}$. For a certain observation $x \in \R^d$ in the data set, such a deviation is usually known as the residual (terminology borrowed from the Statistical Regression literature). In a general framework,  for a given point $x \in \R^d$, we define the \textit{residual} of a  model as a mapping $\varepsilon_x: \R^{d+1} \rightarrow \R_+$, that maps any set of coefficients $\bbeta=(\beta_0, \ldots, \beta_d) \in \R^{d+1}$, into a measure $\varepsilon_x(\bbeta)$ that represents the deviation of the given point $x$ from the hyperplane with those parameters. The larger  this measure, the worse the fitting for such a point $x$. The final goal of fitting an hyperplane for a given set of points $\{x_1, \ldots, x_n\} \subseteq \R^d$ is to find the coefficients minimizing a globalizing function, $\Phi: \R^n \rightarrow \R$, of the residuals of all the points. Equivalently, the fitting problem consists in locating a hyperplane minimizing the globalizing function $\Phi$ of the distances from the demand points to the hyperplane. Different choices for the residuals and the globalizing criteria  will give, in general, different optimal values for the parameters and thus different properties for the resulting hyperplanes. This problem is not new and some of these fitting criteria, as the minisum, minimax and some other robust versions, have been analyzed from a Locational analysis perspective  (see \cite{carrizosa-plastria95,MT83,schobel1,schobel2,schobel3,schobelbook}, among other).

The most natural approach to locate a hyperplane is to consider that residuals, with respect to given points, are individual measures of error and thus,  each residual should be minimized independently of the remaining. Obviously, this approach gives rise to a multicriteria problem \cite{CCFMP95,NW2007}. It is clear that this simultaneous minimization will not be possible in most of the cases and then several  strategies can be followed: one can try to find the set of Pareto fitting curves \cite{CCFMP95} or alternatively, to apply an aggregation function that incorporates the holistic  preference of the Decision-Maker on the different residuals. This last choice is very difficult and the usual approach is to apply the principle of complete uncertainty leading to additive aggregations.

The most popular methods to compute the coefficients of an optimal hyperplane consider that the residuals are the differences from one of the coordinates of the space (which are usually known as vertical/horizontal distances). In this paper we present a new framework for optimally locating/fitting hyperplanes to a set of points that allows the decision-maker to decide within a wide family of residuals and criteria which is the ``best'' for a given sample of data. One of the main contributions of our proposal is the use of modern mathematical programming tools to solve the problems which are involved in the computation of the parameters of the fitting models. The  optimization  models for those problems range from continuous convex programming (CP) to mixed integer nonlinear programming (MINLP) through linear programming (LP). Many of the formulations described in this paper have been implemented in \r{} in order to be available for data analysts.

The framework in this paper introduces a  family of combinations residuals-criteria that allows a great flexibility to accommodate hyperplanes to set of points \cite{NP05,Marin2009}. This new framework can be easily combined with  some of  the mathematical programming techniques for feature selection, to ``choose'' a fixed number of coordinates to explain the dependence between the different dimensions \cite{BM14}, with classification schemes  \cite{BS07}, or  when the coefficients of the linear manifold are required to fulfill a set of linear equations/inequalities.  This framework can also  accommodate general forms of regularization, as upper bound on the $\ell_2$-norm of the coefficients \cite{HK88}, since it would only mean to add additional constraints to the mathematical programming formulations proposed in the paper. The complexity of solving the resulting model depending on the difficulty of the considered regularization constraints.

In order to compare the \textit{goodness of the fitting}  for the different models we have developed a  new generalized measure of fit. This task becomes difficult when one tries to compare  fitting hyperplanes which are built based on different paradigms and purposes. The new measure is provided in order to make meaningful comparisons. This proposal is based on a generalization of the classical coefficient of determination, that will allow to  measure how good is  an optimal hyperplane with respect to the best constant model, $X_d=\beta_0$. This measure will extend the standard coefficient of determination for least squares fitting. We also perform an extensive series of experiments to validate the application of our results applied with different objectives to several set of data.

In our framework, errors are measured as shortest distances, based on a norm, between the given points and the fitting surface. This makes the location problem geometrically invariant which is an interesting advance with respect to vertical/horizontal residuals. Through the paper we observe that this framework also subsumes as particular cases the  standard location methods that consider residuals based on vertical distances (commonly used in Statistics); as well as most of the particular cases of fitting linear bodies using vertical distances but different aggregation criteria described  in the literature, as $\ell_p$ fitting ($\ell_p$-norm criterion),  least quantile of squares \cite{rousseeuw84,BM14}, least trimmed sum of
squares \cite{rousseeuw83,atkinson99},  etc.   As previously mentioned, the problem of optimally locating an hyperplane with respect to a set of demand points is closely related to the estimation phase in multivariate linear regression, where several methods have already been proposed. However,  the use of nonstandard residuals is not usual in the literature of regression analysis although  orthogonal ($\ell_2$) residuals have been already used, see e.g. Euclidean Fitting \cite{BH93,cavalier91,pinsonetal} or Total Least Squares \cite{tls}, mainly applied to bidimensional data. Quoting the reasons for that fact given by Giloni and Padberg in \cite{giloni-padberg02}: ``we have left out a summary of linear regression models using the more general $\ell_\tau$,-norms with $\tau \not\in \{1,2, \infty\}$ for which the computational requirements are considerably more burdensome than in the linear programming case (as they generally require methods from convex programming where machine computations are far more limited today).''

The paper is organized as follows. In Section \ref{sec:1} we introduce the new framework for fitting hyperplanes as well as some results that allows to interpret the results for practical purposes. Next, a residual-aggregation dependent  goodness of fitting index is defined and it is presented an efficient approach for its computation. Section \ref{sec:2} is devoted to the analysis of the classical location methods under the new framework, more precisely, mathematical programming models for adequate aggregation criteria and residuals are provided for: 1) least sum of squares; 2) least absolute deviation; 3)  least quantile of squares and 4) least trimmed of squares fitting. In Sections \ref{sec:3} and \ref{sec:4} we present new methods for the location of hyperplanes assuming that the residuals are  measured as the smallest norm-based distance between the given points (data set) and the linear fitting body using polyhedral norms (Section \ref{sec:3}) and $\ell_{\tau}$ norms (Section \ref{sec:4}), respectively. We also present, in Section \ref{sec:4}, outer an inner approximations for solving the resulting  MINLP problems for $\ell_p$-norms residuals. Finally, Section \ref{sec:5} is devoted to the computational experiments. We report results for synthetic data and for the classical data set given in \cite{dw}.

\section{A flexible methodology for the location of hyperplanes}
\label{sec:1}

Given is a set of  of $n$ observations or demand points (depending that we use the \textit{jergon} of data analysis or location analysis, respectively)  in a $(d+1)$-dimensional space, $\{x_1, \ldots, x_n\} \subset \R^{d+1}$ (we will assume,  for a clearer description of the models, that the first, the $0-th$, component of $x_i$ is the one that account for the intercept in the model,  being $x_{10}=\cdots=x_{n0}=1$). Next, we analyze ways of fitting these observations to a linear form (hyperplane).   For any $y\in \R^{d+1}$,  we shall denote $y_{-0}=(y_1,\ldots,y_d)$, i.e. the vector with the last $d$ coordinates of $y$ excluding the first one. We consider here a flexible framework for the problem of locating/fitting hyperplanes that includes as special cases the classical and most modern   models found in the specialized literature. First, we assume that the point-to-hyperplane deviation is modelled  by a residual mapping  $\varepsilon_x: \R^{d+1} \rightarrow \R_+$, $\varepsilon_x(\bbeta)=\d(x_{-0},\mathcal{H}(\bbeta))$, being $\d$ a distance measure in $\R^d$. This residual represents how ``far'' is the point (observation) $x \in \R^{d+1}$ with respect to the hyperplane $\mathcal{H}(\bbeta) = \{y \in \R^d: (1,y^t) \bbeta = 0\}$ (Some times we will write the hyperplane as $\bbeta^t X = 0$, with   $\bbeta = (\beta_0, \beta_1, \ldots, \beta_d)^t \in \R^{d+1}$.)

%$\varepsilon_x(\bbeta)=\d(x_{-0},\mathcal{H}(\bbeta))$, being $\d$ a distance measure in $\R^d$. Hence, the linear model to be fitted can be written as:
%$$
%\bbeta^t X = 0, \quad \mbox{ with }  \bbeta = (\beta_0, \beta_1, \ldots, \beta_d)^t \in \R^{d+1}.
%$$

Furthermore, the residuals for  each  demand point are aggregated using a globalizing function  $\Phi: \R^n \rightarrow \R$, which for a set of residuals $\bepsilon_1, \ldots, \bepsilon_n$ gives an overall measure of the deviations of  the whole data set with respect to the hyperplane. With this setting, ones  tries to minimize such a globalizing measure of the residuals with respect to all the given demand points.

With this notation, the \textit{Fitting Hyperplane Problem} (FHP) consists in finding $\hat \bbeta\in \R^{d+1}$ such that:
\begin{equation}
\label{eq1}\tag{${\rm FHP}(\Phi, \bepsilon)$}
\hat \bbeta \in \mbox{arg}\min_{\bbeta \in \R^{d+1}} \Phi(\bepsilon_x(\bbeta)),
\end{equation}
where $\bepsilon_x(\bbeta) = (\varepsilon_{x_1}(\bbeta), \ldots, \varepsilon_{x_n}(\bbeta))^t$ is the vector of residuals.

Note that the difficulty of solving \ref{eq1} depends of the expressions for the residuals and the aggregation criterion $\Phi$. If $\Phi$ and $\varepsilon_x$ are linear, the above problem becomes a linear programming problem. In this paper, we consider  a general family of aggregation criteria that includes as particular cases most of the classical ones  used in the literature. Some of  those criteria have  been already considered  for the sake of outlier detection \cite{phonecal,yager2010} or as robust  alternatives to the standard linear regression approach \cite{BM14,giloni-padberg02}.

Let $\lambda_1, \ldots, \lambda_n \in \R$ and let  $\bepsilon \in \R^n$ be the vector of residuals  of all of the demand points in the given data set. We consider aggregation criteria $\Phi: \R^n \rightarrow \R_+$ defined as:

\begin{equation}
\label{phi}
\Phi(\bepsilon) =  \dsum_{i=1}^n \lambda_i\; \bepsilon_{(i)}^p
\end{equation}
where $\bepsilon_{(i)} \in \{\bepsilon_1, \ldots, \bepsilon_n\}$ is such that $\bepsilon_{(1)} \leq \cdots \leq \bepsilon_{(n)}$. Observe that this operator defines a multiparametric family (called \textit{ordered median functions} \cite{NP05}) that depending on the choice of the  $\lambda$-weights captures many of the models proposed in the literature.

Note that the above shape of $\Phi$ is  symmetric and,  for non negative lambda weights, a  monotone function that ensures that the  ordering of the individual residuals do not affect the overall goodness of the fitting. Moreover, it also implies that a componentwise smaller vector of residuals gives rise to a more accurate fitting.

The natural implication of the assumption made about the definition of residuals is that, as expected, the response (projection) of a demand points on a given hyperplane differs from the  classical evaluation and it must be the closest point, with respect to the distance $\d$,  in the located hyperplane $\H(\bbeta)$.
\begin{lem} \label{le:response}
For a given  point $z^t= (1, z_{1}, \ldots, z_{d})$ and the hyperplane $\H(\bbeta)$ the response $\hat z$ consistent with the residual $\bepsilon_{z}=\min_{y\in \mathcal{H}(\bbeta)}\|z_{-0}-y\|$ is given by
\begin{equation*}% \label{eq:normproj}
\hat z= z_{-0} -\frac{\bbeta^t  z}{\|\bbeta_{-0}\|^*}\k(\bbeta),
\end{equation*}
where $\|\cdot \|^*$ is the dual norm to $\| \cdot \|$ and $\k(\bbeta)= \displaystyle \mbox{\rm arg }\dmax _{\|x\|=1} \bbeta_{-0} ^t x$. Moreover,
\begin{equation} \label{eq:normdist}
\bepsilon_{z}=\frac{|\bbeta^t z|}{\|\bbeta_{-0} \|^*}.
    \end{equation}
\end{lem}
\proof
The proof follows applying \cite[Theorem 2.1]{mangasarian} to the definition of residual  $\bepsilon _{z}=\min_{y\in \mathcal{H}(\bbeta)}\|z_{-0}-y\|$.
\endproof

From the above result, the response for a point  with  a  unknown coordinate (w.l.o.g, the last component, $d$),  namely $z=(1,z_{1}, \ldots, z_{d-1},0)^t$, will be given by:
 $$ \hat z_d=  -\frac{\bbeta^t z}{\|\bbeta_{-0}\|^*} \k(\bbeta)_d.
$$
Hence, differentiating $\hat z$ with respect to each $z_j$, $j=1, \ldots, d-1$, we get
$$
\frac{\partial \hat z_d}{\partial z_j} = -\dfrac{\beta_j}{\|\bbeta_{-0}\|^*} \k(\bbeta)_d,
$$
which may be interpreted as the marginal variation of the $d$-th coordinate with respect to  $j$-th  coordinate  whenever the other  dimensions remain constant.

Explicit expressions for such projections, namely, $\ell_1, \ell_\infty$ and $\ell_\tau$-norms, for $\tau>1$  are described in the following lemma.

\begin{lem}
\label{le:responsesnorms}
Let $z= (1, z_{1}, \ldots, z_{d})^t$, then
\begin{enumerate}
\item If $\d$ is the $\ell_1$- distance,
$$
\hat z_k = \left\{\begin{array}{cl} z_{k} & \mbox{ if $|\beta_k| \neq \max \{|\beta_j|: j=1, \ldots, d\}$,}\\
z_k - \frac{\bbeta^t  z}{\|\bbeta_{-0}\|_\infty}v_k, & \mbox{ if $\beta_k = \max \{|\beta_j|: j=1, \ldots, d\}$,}\\
z_k + \frac{\bbeta^t  z}{\|\bbeta_{-0}\|_\infty}v_k, & \mbox{ if $\beta_k = -\max \{|\beta_j|: j=1, \ldots, d\}$,}
\end{array}\right.
$$
for  $k=1, \ldots, d$, and
for  some $v_1, \ldots, v_{d} \geq 0$ such that $\dsum_{j} v_j =1$.
\item If $\d$ is the $\ell_\infty$- distance,
$$
\hat z_k = \left\{\begin{array}{cl}
z_k - \frac{\bbeta^t  z}{\|\bbeta_{-0}\|_1}, & \mbox{ if $\beta_k >0$,}\\
z_k + \frac{\bbeta^t  z}{\|\bbeta_{-0}\|_1}, & \mbox{ if $\beta_k <0$,}
\end{array}\right. \quad k=1, \ldots, d.
$$
%\item If $\d$ is the $\ell_2$- distance,
%$$
%\hat z_k= z_k -\frac{\bbeta^t  z}{\|\bbeta_{-0}\|_2^2} \bbeta_k, \quad k=1, \ldots, d
%$$
\item  If $\d$ is the $\ell_{\tau}$- distance with $1<\tau<+\infty$ then
\begin{equation*} %\label{eq:lpnorm}
\hat z_k = z_{k} -\frac{\bbeta^t  z}{\|\bbeta_{-0}\|_\nu}\k_\tau(\bbeta)_k,\quad k=1, \ldots, d
\end{equation*}
and
$$\k_\tau(\bbeta)_k = \left\{\begin{array}{cl} \frac{\sign(\bbeta_k)   |\bbeta_k|^{\nu/\tau}}{(\sum_{j=1}^d |\bbeta_j|^\nu)^{1/\tau}} & \mbox{if } \bbeta_k\neq 0\\
0 & \mbox{if } \bbeta_k= 0,
\end{array} \right. \quad k=1, \ldots, d,
$$
being $\nu$ such that $\frac{1}{\tau}+\frac{1}{\nu} =1$.
\end{enumerate}

\end{lem}

\proof
The proof  of items 1. and 2.  can be found in \cite{mangasarian}. The proof of item 3. follows from the Lagrangian optimality condition applied to $\dmax_{\|z\|_\tau=1} \bbeta_{-0}\, z$. First, we observe that a Lagrange multiplier exists since the problem is regular at any point of the $\ell_\tau$ unit ball. Next, the Lagrangian function is $L(z,\lambda)=\bbeta_{-0}\, z-\lambda \sum_{k=1}^{d} |z_k|^{\tau}$. Therefore, its partial derivatives are: $\frac{\partial L}{\partial z_k}=\beta_k-\lambda \tau |z_k|^{\tau-1} \sign(z_k)$, for all $k=1,\ldots,d$. Hence, equating to zero the partial derivative, it follows that for any index $k$ such that $z_k^*\neq 0$
\begin{equation} \label{eq:mul-lambda}
\lambda^*=\frac{\beta_k}{\tau|z_k^*|^{\tau-1}} \sign(z_k^*).
\end{equation}

Let us define the sets $I=\{k:\beta_k>0\}$, $J=\{k:\beta_k<0\}$, $K=\{k:\beta_k=0\}$. Now from equation \eqref{eq:mul-lambda}, and taking into account that $\|z\|_\tau=1$, we obtain:
$$ |z_k^*|^\tau=\left\{
\begin{array}{ll} \frac{\left(\sign(z_k^*) \beta_k\right)^\nu}{(\sum_{j=1}^{d} \sign(z_j^*) \beta_j)^\nu} & \mbox{ if } k\in I\cup J \\
0 & \mbox{otherwise.} \end{array} \right.
$$
Moreover, the hessian of $L$ is diagonal and all its entries are negative, namely $\frac{\partial^2 L}{\partial z_k^2}=-\lambda \tau (\tau-1) |z_k^*|^{\tau-2}.$ This implies that $z^*$ and $\lambda^*$ are local maxima.

 In the particular case of $\tau=2$ then one can check that $k_2(\bbeta)_k=\beta_k$ which simplifies the above expression.

\endproof

We note in passing that $\bepsilon _{x}=\d_{\| \cdot \|}(x_{-0},\mathcal{H}(\bbeta))$ and thus, according to the Lemma \ref{le:response}
\begin{equation}
\label{eq:manga}
\d_{\| \cdot \|}(x_{-0}, \H) = \dfrac{|\bbeta^t x|}{\|\bbeta_{-0}\|^*}.
\end{equation}

Observe also that when the demand points in the data set lie exactly on the hyperplane $\mathcal{H}$  all the proposed methods \ref{eq1}  determine the same hyperplane $\mathcal{H}$ as an optimal fitting, for any norm-based residuals while using vertical distance residuals will never produce hyperplanes in the form $\mathcal{H} = \{z \in \R^d: \beta_0 + \beta_1z_1 + \cdots + \beta_{d-1}z_{d-1}=0\}$ since the ``traditional'' methods do not allow zero coefficients for the \textit{dependent}  coordinate. Note also that the vertical distance based methods assume that errors are present only in one of the components (the so-called dependent), so the rest of the variables should be  error-free. In the proposed general framework, this is no longer assumed since there is no distinction between dependent and independent variables for the location/fitting procedure, so errors may be considered in all the components of the points in the given data set.

Remark that the standard residual (vertical distance) is a distance measure that is not induced by a norm, but its expression can be written in a analogous form and so it fits to the shape of the distances that are considered in this paper. In particular, the vertical distance  (with respect to the last coordinate) may be defined as:
\begin{equation}
\label{eq:manga1}
\d_V(x, H) = \dfrac{\left|\beta_dx_d - \dsum_{i=1}^{d-1} \beta_i x_i - \beta_0\right|}{|\beta_d|}.
\end{equation}
\bigskip

The above aggregation criteria (\ref{phi}) and residual functions (\ref{eq:normdist}) are rather general and exhibit good structural properties. On the one hand, they accommodate most of the already considered fitting methods in the literature. On the other hand, one can always exploit its properties and different representations in order to solve the optimization problem \ref{eq1}. In the following we prove some structural properties that imply some sources of solvability of the problem on hands.

For the sake of completeness, we recall the concept of \textit{difference of convex} (D.C.) function. A function $f: \R^d \rightarrow \R$ is said to be a D.C. if there exist $g, h: \R^d \rightarrow \R$ convex functions such that $f$ can be decomposed as the difference between $g$ and $h$. Optimization problems where the objective function and/or the constraints are defined by D.C. functions are called D.C. programming problems and they play an important role in nonconvex optimization because of its theoretical aspects as well as its wide range of applications (see \cite{Thoai99}).

\begin{lem} \label{le:dc}
The globalizing function $\Phi(\bepsilon_x(\bbeta))$ is a D.C. function.
\end{lem}
\proof
In order to prove that the function $\Phi$ is D.C. we will find a convenient representation where we can apply properties of the algebra of D.C. functions. To this for, we introduce the functions:
$$ \varphi _r(\bbeta):= \min \Big\{ \max \{ \bepsilon_{x_{i_1}}(\bbeta)^p,\ldots,\bepsilon_{x_{i_r}}(\bbeta)^p : i_1<i_2<\ldots<i_r,\; \forall i_1,i_2,\dots,i_r\}\Big\},$$
for $r = 1, \ldots, n$, where $\bepsilon_{x}(\bbeta)=\d_{\| \cdot \|}(x_{-0}, \H)$.

It is a simple observation that $\varphi _r(\bbeta)$ coincides with the $p$-power  $r$-th residual sorted in non-decreasing sequence, namely $\varphi _r(\bbeta)=\bepsilon_{x_{(r)}}(\bbeta)^p$ for all $r=1,\ldots,n$. Hence, we get that $\Phi(\bepsilon_x(\bbeta))=\sum_{i=1}^n \lambda_i \varphi _r(\bbeta)$.

To finish the proof it suffices to prove that each function $\varphi_r$ is D.C. since linear combinations of D.C. functions are D.C.. Next, we start analyzing the residual function $\bepsilon_x(\bbeta)=d(x,\mathcal{H}(\bbeta))$. Assuming that $d$ is a norm based distance given in the form of (\ref{eq:manga}) or (\ref{eq:manga1}), one can use those expressions to conclude that for each observation $x$, $\bepsilon_x(\bbeta)$ is D.C. function of $\bbeta$. Raising to the power $p$ with $p\ge 1$ is also D.C., since it is the result of composing with a convex function (observe that residuals are non-negative). Finally, the operations of taking maxima and minima of D.C. functions are closed within this family \cite{Thoai99}. This proves that $\varphi_r$ is D.C. for all $r$ and this concludes the proof.
\endproof

We note in passing that the D.C. character of our globalizing criterion allows the application of all the available results on the optimization of this class of functions (see e.g. \cite{Thoai99}). In spite of that, we can give more efficient representations that may help latter in the resolution of particular hyperplanes. These representations are based on simpler functions which replace $\varphi$ by more friendly classes of functions (with regards to the optimization phase).

\begin{prop}
The globalizing function $\Phi(\bepsilon_x(\bbeta)):=\sum_{i=1}^n \bepsilon_{x_i}(\bbeta)^p+\sum_{r=2}^n (\lambda_r-\lambda_{r-1}) \theta_r(\bbeta),$
where $\theta_r(\bbeta)=\max \Big\{ \bepsilon_{x_{i_1}}(\bbeta)^p+\ldots+\bepsilon_{x_{i_r}}(\bbeta)^p:\begin{array}{c} \{i_1,\ldots,i_r\}\subset \{1,\ldots,n\} \\ i_1<i_2<\ldots<i_r \end{array}\Big\}, \; r=2,\ldots,n.$
(The reader may observe that the functions $\theta_r$ are usually called $r-centrum$ in the specialized literature of optimization (\cite{NP05}).)
\end{prop}
\proof
This representation follows from the combination of the result in Lemma \ref{le:dc} and   \cite[Theorem 3.6]{GNPU11}.
\endproof

The following result states a mathematical programming formulation for the generalized fitting hyperplane problem, for any choice of $\Phi$ and $\bepsilon_x$.

\begin{theo}
Let $\{x_1, \ldots, x_n\} \subseteq \R^{d+1}$ be a  set of demand points, $\lambda \in \R^n_+$,  $p = \dfrac{r}{s}\in \Q$ and $\|\cdot\|$ a norm in $\R^d$. The Problem \ref{eq1} is equivalent to the following mathematical programming problem:

\begin{align}
\min \dsum_{j=1}^n \lambda_j \theta_j\label{eq:form1}\tag{${\rm LR}_{\Phi,\|\cdot\|}$}\\
\mbox{ s.t. } & \bepsilon_i \geq \frac{|\bbeta^t x_i|}{\|\bbeta_{-0}\|^*}, &\forall i=1, \ldots, n,\label{nlc}\\
& z_i \leq \theta_j + M(1-w_{ij}), &\forall i, j =1, \ldots, n,\label{ctr:1}\\
& z_i^s \geq \bepsilon_i^r, & \forall i=1, \ldots, n,\label{ctr:2}\\
& \dsum_{i=1}^n w_{ij} = 1, & \forall j=1, \ldots, n,\label{ctr:3}\\
& \dsum_{j=1}^n w_{ij} = 1, & \forall i=1, \ldots, n,\label{ctr:4}\\
& \theta_j \geq \theta_{j-1}, &\forall j=2, \ldots, n,\label{ctr:5}\\
& w_{ij} \in \{0,1\}, &\forall i=1, \ldots, n,\nonumber\\
&  z, \theta \in \R^n_+, \bbeta \in \R^{d+1}\nonumber.
\end{align}

\end{theo}

Note that the above problem is a mixed integer non linear programming problem, whose continuous relaxation is in general non convex due to the constraints \ref{nlc}.  Apart from the mathematical programming formulation above, one may use alternative (in some cases better) formulations for the ordering problems as those provided in \cite{FPP2014}.
 In particular, some important special ordered median aggregation criteria allow to have a simpler formulation that avoids the use of binary variables.  The following result shows a better formulation for the fitting problem under the assumption that $0\le \lambda_1\le \ldots \le \lambda_n$. We call this setting for lambda the  \textit{monotone case}.

\begin{theo}
%\label{theo:convex}
Let $\{x_1, \ldots, x_n\} \subset \R^{d+1}$ be a  set of demand points, $\lambda \in \R^n$, such that $0\le \lambda_1 \leq \cdots \leq \lambda_n$,  $p = \dfrac{r}{s}\in \Q$ with $r > s \in \N$, $\gcd(r,s)=1$ and $\|\cdot\|$ a norm in $\R^d$. Then, \ref{eq1} is equivalent to the following mathematical programming problem:

\begin{align}
\min \dsum_{j=1}^n v_j+ \dsum_{i=1}^n w_i \nonumber\\%\label{eq:form2}\tag{${\rm LR}^C_{p,\|\cdot\|}$}\\
\mbox{ s.t. } & \eqref{nlc}, \eqref{ctr:2},\nonumber\\
& v_j + w_i \geq \lambda_i z_j, \forall i, j =1, \ldots, n,\nonumber\\%\label{convex:1}\\
& z_i, \theta_i \geq 0, v, w  \in \R^n,  \bbeta \in \R^{d+1}\nonumber.
\end{align}

\end{theo}

\proof
The proof follows by the representation of the ordering between the residuals by permutation variables, which for $0\le \lambda_1 \leq \cdots \leq \lambda_n$, allows to write the objective function in \ref{eq1} as an assignment problem which is totally unimodular, so it can be equivalently rewritten using its dual problem. The interested reader is refereed to \cite{BPE14} for further details on this transformation.
\endproof

The reader may observe that, based on an alternative representation, the nonlinear constraints $z_i^s \geq \bepsilon_i^r$ for all $i=1, \ldots, n$  can be transformed into a set of second order cone constraints using the following result which is a simplified version of Lemma 1 in \cite{BPE14}. This implies that those constraints can be efficiently handled by nowadays nonlinear solvers since they are convex and friendly for the optimization.

\begin{lem} \label{le:n2}
Let $r,s\in\N\setminus\{0\}$ with $\gcd(r,s)=1$, and $k=\lfloor\log_2(r)\rfloor$. Then, there exist variables $u_1, \ldots, u_{k-1} \geq 0$ such that each constraint $z^s \geq \bepsilon^r$ in \ref{eq:form1}  can be equivalently written as constraints in the form:

      $$
      \begin{array}{ll}
        u_{j}^{2} &\leq u_{l}^{a_j} z^{b_j} \bepsilon^{c_j},\\
        \bepsilon^2 & \leq u_h u_{h-1}^{d_h} z^{f_h}\bepsilon^{g_h},\\
        u_j &\geq 0
      \end{array}
      $$
      with $j=1, \ldots, k-1$ and such that $1 \leq a_j+b_j+c_j \leq 2$ for given $a_j, b_j, c_j \in \Z_+$ and $d_h, f_h, g_h \in \Z_+$ such that
$d_h+b_h+c_h =1$.
\end{lem}

By the above lemma, the nonlinear constraints in the form $z^s \geq \bepsilon^r$ are written as second order cone  constraints in the form $X^2 \leq Y Z$ or $X^2 \leq Y$ (for some choices of the variables $X$, $Y$ and $Z$ in our model). These constraints are then equivalent to one of the following two  semidefinite constraints:
$$ \small
\left(\begin{array}{ccc} Y+Z  & 0 & 2X \\ 0 & Y+Z & Y-Z \\ 2X & Y-Z & Y+Z \end{array} \right) \succeq 0, \; Y+Z\ge 0
\mbox{ or  }
\left(\begin{array}{ccc} Y  & 0 & 2X \\ 0 & Y & Y \\ 2X & Y & Y \end{array} \right) \succeq 0, \; Y\ge 0.
$$
Hence, the difficulty of solving Problem  \ref{eq:form1},  depends essentially on the choice of the residuals since all except  constraints \eqref{nlc} are linear or second order cone constraints which can be efficiently handled with nowadays modern optimization techniques. In the next sections we analyze different choices of the residuals.

\begin{rem}[Subset Selection and Regularization]
In the case where the number of points ($n$) is much smaller than the dimension of the space ($d$), it is common in Statistics to compute fitting hyperplanes over a smaller dimensional space. The new space is determined by those components that, after projecting, allows a good fitting when it is compared to the dimension of the new space. Several methods have been proposed in the recent literature to perform such a computation. If the dimension of the new space, $q < d$, is given, a constraint in the form $\|\beta_{-0}\|_0\leq q$ (here $\|\cdot\|_0$ stands for the support function or nuclear norm, i.e., the number of nonzero components of the vector) may be included in the mathematical programming formulation (see \cite{miller,Bertsimas2016}), which gives rise to the so called Subset Selection Problem. If such a dimension is not known, regularization methods that penalize the number of nonzero elements or the size of $\beta_{-0}$ can be applied to solve the Feature Selection Problem (see \cite{miyashiro15}). Note that both types of approaches can be easily incorporated in our models.
\end{rem}

\subsection{Goodness of Fitting}
After addressing the problem of locating/fitting a hyperplane with respect to a set of points,  we will analyze the goodness of this fitting extending the well known coefficient of determination in Regression Analysis. For the sake of presentation, we assume that the variable that needs to be analyzed in terms of dependence to the others is the last coordinate $X_d$, or in other words $Y=X_d$. The  \textit{goodness of fitting index} is defined as:
\begin{equation*}%\label{gcod}
\GCoD_{\Phi, \bepsilon} = 1- \dfrac{\Phi^* }{\Phi_0^*},
\end{equation*}
where $\Phi^*$ is the optimal value of \eqref{eq1}, namely $\Phi(\bepsilon_x(\hat \bbeta))$, and $\Phi^*_0$ is the optimal value of \ref{eq1} when it is additionally required that $\bbeta$ is in the form $\bbeta = (\beta_0, \overbrace{0, \ldots, 0}^{d-1}, -1)$,  i.e.  the hyperplane is imposed to be constant ($X_d = \beta_0$). Note that the components $1, \ldots, d-1$ do not appear in the model. Hence, $\Phi^*_0$ measures the global error assumed by the \textit{best} fitting  ``vertical'' hyperplane; %$\beta_0 + \beta_d X_d = 0$
whereas  $\GCoD_{\Phi, \bepsilon}$ measures the improvement of the model that considers all the dimensions with respect to the one that omits all  (except one) of  them . Observe that this coefficient coincides with the classical coefficient of determination provided that the aggregation criteria is the overall sum and the residuals are the squared vertical distances: in that case $\widehat{\beta}_0= \overline{x}_{\cdot d}$ (the sample mean of the \textit{dependent} variable).

The $\GCoD$ clearly verifies one of the important properties of the standard coefficient of determination, $0 \leq \GCoD_{\Phi, \bepsilon} \leq 1$. Furthermore, one may interpret the coefficient as a measure of how good is the best possible  hyperplane under certain criterion and residual choice with respect to the best \textit{horizontal} hyperplane. When $\GCoD$ is close to $0$, it is because $\Phi^* \simeq \Phi_0^*$, so not appreciable improvement is given by the complete model (which considers all the components) with respect to the simple constant model; whenever $\GCoD$ is close to $1$, it means that $\Phi^* \ll \Phi_0^*$, being the proposed model significatively better than the constant model (note that $\GCoD=1$ iff $\Phi^*=0$, i.e., when the model perfectly fits the demand points). Hence, the closer the $\GCoD$  to one, the better the fitting; whereas the  closer to zero,  the better is the constant model with respect to the full model.

Observe that the above definition coincides with some of the choices to measure the goodness of fitting for robust alternatives to the least sum of squares methodology (see \cite{MS87}).

To obtain the $\GCoD$, apart from solving \ref{eq1} to get $\Phi^*$,  we must also solve the problem:
\begin{equation}
\label{eq2-0}
\Phi_0^* =  \min_{\beta_0 \in \R} \Phi(\d(x_1,\H_0), \ldots, \d(x_n,\H_0)),
\end{equation}
 where $\H_0 = \{y \in \R^d: y_d = \beta_0\}$ for some $\beta_0 \in \R$.

\begin{lem} \label{le:ompres}
If the  residual mapping $\varepsilon_x: \R^{d+1} \rightarrow \R_+$ is  induced by a norm $\|\cdot\|$. Then,
Problem (\ref{eq2-0}) is equivalent to
\begin{equation}
\label{eq2}\tag{${\rm LRP}^0_{\Phi, \bepsilon}$}
\Phi_0^* =  \min_{\beta_0 \in \R} \Phi(\kappa_{\bepsilon} |x_{1d}-\beta_0|, \ldots, \kappa_{\bepsilon} |y_{nd}-\beta_0|),
\end{equation}
where
\begin{equation*}%\label{eq:kappa}
\kappa_\varepsilon = \dfrac{1}{\dmax_{z \in \R^d: \|z\|\leq 1} z_d}
\end{equation*}
\end{lem}
\proof
For the point $x_k$ in the data set, the residual under the assumption $X_d=\bbeta_0$ is $\varepsilon_{x_k} (\beta_0) = \d(x_k, \mathcal{H}_0) = \min_{y \in \mathcal{H}_0} \|x_k-y\|$, where $\mathcal{H}_0 = \{y \in \R^d: y_d = \beta_0\}$ for some $\beta_0 \in \R$. Then, by (\ref{eq:normdist}) in Lemma \ref{le:response}
$$
\varepsilon_{x_k} (\beta_0) = \dfrac{|x_{kd} - \beta_0|}{\|(0, \ldots, 0, -1)\|^*}
$$
with $\|\cdot\|^*$ the dual norm of $\|\cdot\|$. By definition of the dual norm $\|y\|^* = \dmax_{z \in \R^d: \|z\|\leq 1} z^t y$. Hence, applying such a definition to $y=(0, \ldots, 0, -1)$ the result follows.
\endproof

From the above result it is easy to see that $\kappa_\varepsilon=1$, provided that $\varepsilon_x$ is induced by any $\ell_p$ norm, even for the $\ell_1$ and the $\ell_\infty$ cases. However, as we will see in Section \ref{sec:3}, not all the norms have the same $\kappa_\varepsilon$ constant.

\medskip

Next, with our specifications for $\Phi$, given by \ref{eq1},  the problem to be solved  to obtain $\Phi^*_0$ is:
\begin{equation}
\label{eq3}\tag{${\rm LRP}^0_{\lambda, p}$}
\Phi_0^* = \kappa_{\bepsilon} \min_{\beta_0 \in \R}  \flp(\beta):= \dsum_{i=1}^n \lambda_i\; \varepsilon_{(i)}^p
\end{equation}
where $\varepsilon_i = |x_{id}-\beta_0|$ for $i=1, \ldots, n$.

Solutions to Problem \ref{eq3} for a given $\beta_0 \in \R$ motivate the introduction of the concept of \textit{ordered median point}. Indeed, $\beta_0$ is a $(\lambda, p)$-\textit{ordered median point} ($(\lambda, p)$-omp in short)  if it is an optimal solution to \ref{eq3}.%\cite{NP99}.

Some special cases of $(\lambda, p)$-omp are well-known and widely used in the so-called location analysis literature. If $\lambda_i=1$ for all $i=1, \ldots, n$, the $(\lambda, 1)$-omp is  known  to coincide with the median, ${\rm median}(x_{1d}, \ldots, x_{nd})$, of $\{x_{1d}, \ldots, x_{nd}\}$;  while the $(\lambda, 2)$-omp is the arithmetic mean of the $x_{.d}$-values.

In the general case, i.e. for arbitrary $\lambda$ and $p$, the ordered median points do not have closed form expressions \cite{FPP2014,FPPS2016}, although they have been around in the field of Location Analysis for several years \cite{NP99,NP05}. Moreover, they can be  obtained, as shown below, to be used in the computation of the goodness of fitting index.

In the following we show how to solve \ref{eq3} for general choices of non-negative vectors $\lambda$ and $p\in [1,+\infty)$. Without loss of generality we assume that $x_{1d}\le x_{2d}\le \ldots \le x_{nd}$. Let us denote further by $\alpha_{ik}:= \frac{x_{id}+x_{kd}}{2}$ the solution of the equation $\bepsilon_i^p(\bbeta)=\bepsilon_ k^p(\bbeta)$ for all $i<k, \; i,k=1,\ldots,n$  in the range $(x_{1d},x_{nd})$. Let $\mathcal{A}$ be the set containing all the $x_{.d}$ and $\alpha$ points and denote by $z_k$ the $k$-th point in $\mathcal{A}$ sorted in non-decreasing sequence.  By construction, in the interval $I_k=(z_k,z_{k+1})$ all the functions $\bepsilon_i^p(\bbeta)$ are monotone for all $i=1,\ldots,n$.

\begin{lem}
The function $\flp(\bbeta)$ has at most one critical point $\bbeta^*\in I_k$.
\end{lem}
\proof
For all $\bbeta\in I_k$, the function $\flp$ is a non-negative linear combination of monotone functions. Therefore, its derivative can vanish in at most one point.
\endproof

Let us denote by $\mathcal{A}_c$ the set of all the critical points of the function $\flp$ in the interval $(x_{1d},x_{nd})$. Observe that the cardinality of this set is $O(n^2)$.

\begin{theo}\label{teor:ompset}
For any non-negative vector $\lambda$ and $p\in [1,\infty)$ the set $ \mathcal{A}\cup \mathcal{A}_c$ always contains a $(\lambda,p)$-omp.
\end{theo}
The reader may observe that the implication of the above theorem is that the $\hat \bbeta_0$ value can be always obtained by a simple enumeration of the set $ \mathcal{A}\cup \mathcal{A}_c$. Then, $\Phi^*_0=\kappa_{\bepsilon}\sum_{i=1}^n\lambda_i|x_{id}-\hat \bbeta_0|^p_{(i)}$. Thus, the complexity of computing $\GCoD$  is essentially the same as the resolution of Problem \ref{eq1}, which must be solved to obtain $\Phi^*$.

\section{Classical Methods under the new framework}
\label{sec:2}

In this section we show how several classical models of fitting with hyperplanes  can be cast into our general framework. We assume that we are given a set of points $\{x_1, \ldots, x_n\} \subseteq \R^{d+1}$. In classical models in the literature, the residuals are defined as the vertical distance (with respect to the last coordinate) from the point to the hyperplane:
\begin{equation}
\label{eq:cla-res}
\bepsilon_x (\bbeta) = \left|x_d - \dsum_{k=0}^{d-1} \frac{\beta_k}{\beta_d} x_k\right|.
\end{equation}
Therefore, the difference between the considered models comes from the choice of the globalizing criterion  $\Phi$ that aggregates the residuals. We have pointed out in the previous section that an important factor, in determining the difficulty of solving the mathematical programming problems for the fitting model, is the choice of the residual. This element influences  much more the difficulty of the problem than the globalizing criterion. We shall show in this section how to handle, within this framework,  the following 4 well-regarded models: Least Sum of Squares (LSS), Least Sum of Absolute Deviation (LAD), Least Quantile of Squares (LQS)  and  Least Trimmed Sum of Squares (LTS). These four well-known models are presented below as particular cases of our general framework described in \ref{eq1}.

A particularity of the models where the residuals are measured as the vertical distance between the point and the hyperplane, is that the response for a given data $z$ coincides with $\hat z = z_d - \sum_{k=0}^{d-1} \frac{\widehat{\beta}_k}{\widehat{\beta}_d} z_k$, which is the direct evaluation of $z$ over the linear function that defines the fitted hyperplane. This property will not be valid, in general, for  residuals different from the vertical distance.

\subsection{Least Sum of Squares fitting problem}

We start our analysis with the LSS method, credited to Gauss. It is the most widely used approach to estimate the coefficients of a linear model because its simplicity  and its theoretical implications for the inference over the total population. However, somehow restricting hypotheses are required in order to be  applied (see e.g. \cite{giloni-padberg02}).

The LSS criterion is defined as the sum of the squares of the residuals, that is:
$$
\Phi_{LSS}(\bepsilon_1, \ldots, \bepsilon_n) = \dsum_{i=1}^n \bepsilon_i^2,
$$
where the residuals $\bepsilon_k$ are given by \eqref{eq:cla-res}.

In case $n>d$,  assuming without loss of generality that $\beta_d=1$, and that the given points are linearly independent, the optimality conditions of the problem allow to compute the best LSS parameters as:
$$
\bbeta = (X^t X)^{-1} X^t y
$$
where $X$ is the $n\times d$-matrix obtained from the sample data by columns and $y^t = (x_{1d}, \ldots, x_{nd}) \in \R^{n}$ are the responses of the  last component of the model. Hence, the complexity of computing the parameters under the LSS method is $O(nd^2)$ which results from the complexity of multiplying $n\times d$ matrices. However, even though there is a closed form formula, it may appear numerical errors when computing the inverse of the matrix $X^t X$ if the rows of $X$ are linearly dependent or close to the linear dependence.
Alternatively, one can compute the parameter $\bbeta$, regardless of the degree of dependence of the variables in the model by solving either a quadratic programming or a second order cone programming problem; which is nowadays doable with on-the-shell software.

\begin{theo}
 An optimal parameter $\bbeta^* \in \R^{d}$ that minimizes $\Phi_{LSS}$ can be obtained by solving any of the following two problems:

\hspace*{-1.5cm}\begin{tabular}{l|l}
\parbox{0.5\linewidth}{%
\begin{align}
&\min \; \dsum_{i=1}^n  z_i^2,\label{eq:olsQP}\tag{${\rm LSS}_{\rm QP}$}\\
&z_i \geq x_{id} - \dsum_{k=1}^{d-1} \beta_k x_{ik} - \beta_0,\nonumber\\% \label{ctr:8}
&\bbeta \in \R^{d}, z \in \R^n_+.\nonumber
\end{align}}
&
\parbox{0.5\linewidth}{%
\begin{align}
&\min\; \dsum_{i=1}^n  w_i, \label{eq:olsSOCP}\tag{${\rm LSS}_{\rm SOCP}$}\\
&w_i \geq \left( x_{id} - \dsum_{k=1}^{d-1} \beta_k x_{ik} - \beta_0\right)^2, \label{ctr:9}\\%
&\bbeta \in \R^{d}, z  \in \R^n_+.\nonumber
\end{align}}
\end{tabular}
\end{theo}

\proof
Denote by $z_i = x_{id} -  \beta_0 - \dsum_{k=1}^{d-1} \beta_k x_{ik}$ and by $w_i = \left( x_{id} - \dsum_{k=1}^{d-1} \beta_k x_{ik} -\beta_0 \right)^2$, for $i=1, \ldots, n$. Note that  the objective functions in \eqref{eq:olsQP} and \eqref{eq:olsSOCP} coincide:
$$
\dsum_{i=1}^n z_i^2 = \dsum_{i=1}^n (x_{id} - \dsum_{k=1}^{d-1} \beta_k x_{ik} - \beta_0)^2 = \dsum_{i=1}^n w_i
$$
Next, the minimization character of the objective function allows us to relax the equality constraint definition of the auxiliary variables to $\geq$-constraints and then the result follows.
\endproof
The reader may observe that LSS corresponds to \ref{eq1} with $\lambda^t=(1,\ldots,1)$, $p=2$ and $\bepsilon$ the vertical distance.
\subsection{Least Absolute Deviation  fitting problem}

Another well-explored choice of residuals and criterion  is the so called LAD method, introduced by Edgeworth in 1887. The globalizing criterion is the sum of the absolute value of the vertical residuals:
 $$
\Phi_{LAD}(\bepsilon_1, \ldots, \bepsilon_n) = \dsum_{i=1}^n |\bepsilon_i|.
$$
Note that LAD corresponds to the model \ref{eq1} for with $\lambda^t=(1,\ldots,1)$ and $p=1$. The optimal coefficients obtained with this method are known to be more robust than those by the LSS method. It follows that the mathematical programming model to be solved under this choice is:

\begin{equation}
\label{eq:lad0}
\min_{\bbeta \in \R^{d+1}} \dsum_{i=1}^n \left|x_{id}  -\beta_0 - \dsum_{k=1}^{d-1} \beta_k x_{ik} \right|
\end{equation}
(assuming w.l.o.g. that $\beta_d=1$).

Observe that the above problem to compute the best LAD hyperplane can be actually formulated as a linear programming problem by replacing in  \eqref{eq:olsSOCP} the quadratic constraints by those which model the absolute value.

\subsection{Least Quantile of Squares  fitting problem}

Next, we describe another method known as Least Quantile of Squares, recently introduced by Bertsimas and Mazumder \cite{BM14}, which is a generalization of the Least Median of Squares (LMS) introduced by Hampel (1975). It also considers vertical distances as residuals, but the residuals are aggregated to minimize the $r$-quantile of the distribution of residuals ($r$ can range in $\{1,\ldots,n\}$).
 $$
 \Phi_{LQS}(\bepsilon_1, \ldots, \bepsilon_n) = r-{\rm quantile}(\bepsilon_1^2, \ldots, \bepsilon_n^2):=  \bepsilon_{(r)}^2.
 $$
which also fits to the general form of the aggregating criteria considered in this paper. In this case, following the notation introduced in \eqref{phi}, the LQS hyperplane can be obtained for $p=2$ and $\lambda=(\overbrace{0, \ldots, 0}^{(r-1)}, 1, \overbrace{0, \ldots, 0}^{(n-r)})$. (Observe that
LMS hyperplane is also obtained within the same scheme when $p=2$ and $\lambda=(\overbrace{0, \ldots, 0}^{\lfloor\frac{n}{2}\rfloor},$ $1, \overbrace{0, \ldots, 0}^{\lfloor\frac{n}{2}\rfloor})$.)

\begin{theo}
An optimal parameter $\bbeta^* \in \R^{d+1}$ for LMS method can be obtained by solving the following problem:
\begin{align}
&\min \mbox{ } \quad \theta_{r} \label{eq:LMS}\tag{${\rm LMS}_{\rm IP}$}\\
\mbox{ s.t. } \mbox{ } \quad & \eqref{ctr:1}, \eqref{ctr:3}-\eqref{ctr:5}, \eqref{ctr:9}, \nonumber\\
&\bbeta \in \R^{d}, z, \theta \in \R_+^n, w_{ij} \in \{0,1\}, \forall i,j=1, \ldots, n,\nonumber
\end{align}
\end{theo}
\subsection{Least Trimmed Sum of Squares  fitting problem}
Finally, we present analogous formulations for the LTS method. This method was introduced by Rousseeuw \cite{rousseeuw84} as a very robust alternative to the LSS method, in that it has a high breakdown point. With our notation, the residuals are again considered as the vertical distance, $p=2$  but the aggregation criterion is now:
$$
\Phi_{LTS}(\bepsilon_1, \ldots, \bepsilon_n) =  \dsum_{i=1}^h \bepsilon_{(i)}^2
$$
where $\bepsilon_{(i)} \in \{\bepsilon_1, \ldots, \bepsilon_n\}$ with $\bepsilon_{(i)} \leq \bepsilon_{(i+1)}$ for $i=1, \ldots, n-1$, and $h \in \{1, \ldots, n\}$. Note that in this problem one tries to minimize the sum of the $h$ smallest squared residuals, discarding the remaining, and then, adjusting the model to the $h$ \textit{closest} points. The most common choice for $h$ is $\lfloor \dfrac{n}{2} \rfloor$, considering the best $50\%$ square residuals to compute the hyperplane (thus, discarding the other $50\%$ of the data). The choice of $h$ allows to control which part of the data set are \textit{sacrificed} to find a better hyperplane. We denote by $LTS(\alpha)$ the LTS method when $100 - \alpha \%$ of the data is discarded, i.e., the percentage of the data that may be considered as outliers.

A suitable mathematical programming formulation for the $LTS(\alpha)$ method is stated in the following result.

\begin{theo}
An optimal parameter $\bbeta^* \in \R^{d+1}$ for $LTS(\alpha)$ method can be obtained by solving  the following problem:
\begin{align}
\min \;& \dsum_{i=1}^{\lceil \alpha n\rceil} \theta_{j} \label{eq:LTS}\tag{${\rm LTS}(\alpha)_{\rm IP}$}\\
\mbox{ s.t. } \mbox{ } \quad & \eqref{ctr:1}, \eqref{ctr:3}-\eqref{ctr:5}, \eqref{ctr:9}, \nonumber\\
&\bbeta \in \R^{d+1}, z, \theta \in \R^n,  w_{ij} \in \{0,1\}, \forall i,j=1, \ldots, n.\nonumber
\end{align}

\end{theo}

We illustrate the differences of the above classical models in a well-known data set that appears in \cite{phonecal}. The algorithms were implemented in \r~ with the \texttt{Gurobi} callable library.

\begin{ex}
\label{ex:stars}
The data considered in this example consists of $47$ points in $\R^2$ about stars of the CYG OB1 cluster in the direction of Cygnus \cite{stars}.  The first coordinate, $X_1$, is the logarithm of the effective temperature at the surface of the star and the second one, $X_2$, is the logarithm of its light intensity. This data set has also been analyzed in \cite{phonecal} and \cite{yager2010}, among others.

We run the LSS, LAD, LMS and LTS($\alpha$) with $\alpha \in \{25, 50, 75, 90\}$. The obtained lines and the goodness of fitting ($\GCoD_{\Phi, \bepsilon}$) are shown in Figure \ref{f:fig3}.

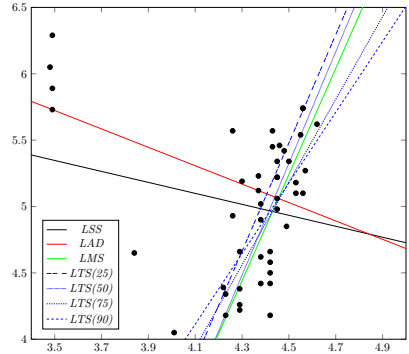
\begin{figure}[t!]
\begin{center}
    \begin{subfigure}[c]{0.5\linewidth}
\scalebox{1}{\begin{tabular}{ccc}
Method & Line & $\GCoD$\\\hline\hline
LSS &  y=-0.4133 x + 6.7934 &  0.0442\\
LAD & y= -0.6931 x + 8.1492 & 0.0065 \\
LMS & y = 4 x -127.6 & 0.0765 \\
LTS(25) & y= 4.0767 x -12.8668 & 0.7328\\
LTS(50) & y=4.2105 x  -13.6231 & 0.6057\\
LTS(75) & y= 3.1176 x  -8.8461 & 0.7702 \\
LTS(90) & y = 2.6620 x -6.8016 & 0.6751\\\hline
\end{tabular}}
\end{subfigure}%
~
\begin{subfigure}[c]{0.5\linewidth}
\scalebox{0.8}{%\begin{figure}[h]
%\centering
\centering
\begin{tikzpicture}[scale=0.55]
  \begin{axis}[
  xmin=3.4,   xmax=5.0,
	ymin=4.0,   ymax=6.5,
xtick = {3.5, 3.7, 3.9, 4.1, 4.3, 4.5, 4.7, 4.9},
ytick = {4.0, 4.5, 5.0, 5.5, 6.0, 6.5},
legend pos=south west,
    %anchor=origin,  % Align the origins
    x=7cm, y=4cm,   % Set the same unit vectors
    ]

        \addplot[mark=none, black, thick] plot coordinates {
        (5.0, -0.4133039*5.0+6.793467)
                (3.4, -0.4133039*3.4+6.793467)
    };

         \addplot[mark=none, red, thick] plot coordinates {
        (8.149205/0.6931818, 0)
                (3.4, -0.6931818*3.4+8.149205)
    };

       \addplot[mark=none, green, thick] plot coordinates {
        (12.76/4.0,0)
                (5.0, 4*5.0-12.76)
    };

     \addplot[mark=none, blue,thick, dash pattern=on 4pt off 1pt on 4pt off 4pt] plot coordinates {
        (12.86685/4.076726, 0)
                (5.0, 5.0*4.076726-12.86685)
    };

     \addplot[mark=none, blue,thick, dash pattern=on 0.5pt off 0.5pt on 0.5pt off 0.5pt] plot coordinates {
        (13.62316/4.210526, 0)
                (5.0, 5.0*4.210526-13.62316)
    };

     \addplot[mark=none, blue,thick, dash pattern=on 1pt off 1pt on 1pt off 1pt] plot coordinates {
        (8.846176/3.117647, 0)
                (5.0, 5.0*3.117647-8.846176)
    };

      \addplot[mark=none, blue,thick, dash pattern=on 2pt off 2pt on 2pt off 2pt] plot coordinates {
        (6.801652/2.662076, 0)
                (5.0, 2.662076*5.0-6.801652)
    };

    \legend{LSS, LAD, LMS, LTS(25), LTS(50), LTS(75), LTS(90)}

\addplot[mark=*, very thin] plot coordinates {
 ( 4.37 , 5.23 )
 };
\addplot[mark=*, very thin] plot coordinates {
 ( 4.56 , 5.74 )
 };
\addplot[mark=*, very thin] plot coordinates {
 ( 4.26 , 4.93 )
 };
\addplot[mark=*, very thin] plot coordinates {
 ( 4.56 , 5.74 )
 };
\addplot[mark=*, very thin] plot coordinates {
 ( 4.3 , 5.19 )
 };
\addplot[mark=*, very thin] plot coordinates {
 ( 4.46 , 5.46 )
 };
\addplot[mark=*, very thin] plot coordinates {
 ( 3.84 , 4.65 )
 };
\addplot[mark=*, very thin] plot coordinates {
 ( 4.57 , 5.27 )
 };
\addplot[mark=*, very thin] plot coordinates {
 ( 4.26 , 5.57 )
 };
\addplot[mark=*, very thin] plot coordinates {
 ( 4.37 , 5.12 )
 };
\addplot[mark=*, very thin] plot coordinates {
 ( 3.49 , 5.73 )
 };
\addplot[mark=*, very thin] plot coordinates {
 ( 4.43 , 5.45 )
 };
\addplot[mark=*, very thin] plot coordinates {
 ( 4.48 , 5.42 )
 };
\addplot[mark=*, very thin] plot coordinates {
 ( 4.01 , 4.05 )
 };
\addplot[mark=*, very thin] plot coordinates {
 ( 4.29 , 4.26 )
 };
\addplot[mark=*, very thin] plot coordinates {
 ( 4.42 , 4.58 )
 };
\addplot[mark=*, very thin] plot coordinates {
 ( 4.23 , 3.94 )
 };
\addplot[mark=*, very thin] plot coordinates {
 ( 4.42 , 4.18 )
 };
\addplot[mark=*, very thin] plot coordinates {
 ( 4.23 , 4.18 )
 };
\addplot[mark=*, very thin] plot coordinates {
 ( 3.49 , 5.89 )
 };
\addplot[mark=*, very thin] plot coordinates {
 ( 4.29 , 4.38 )
 };
\addplot[mark=*, very thin] plot coordinates {
 ( 4.29 , 4.22 )
 };
\addplot[mark=*, very thin] plot coordinates {
 ( 4.42 , 4.42 )
 };
\addplot[mark=*, very thin] plot coordinates {
 ( 4.49 , 4.85 )
 };
\addplot[mark=*, very thin] plot coordinates {
 ( 4.38 , 5.02 )
 };
\addplot[mark=*, very thin] plot coordinates {
 ( 4.42 , 4.66 )
 };
\addplot[mark=*, very thin] plot coordinates {
 ( 4.29 , 4.66 )
 };
\addplot[mark=*, very thin] plot coordinates {
 ( 4.38 , 4.9 )
 };
\addplot[mark=*, very thin] plot coordinates {
 ( 4.22 , 4.39 )
 };
\addplot[mark=*, very thin] plot coordinates {
 ( 3.48 , 6.05 )
 };
\addplot[mark=*, very thin] plot coordinates {
 ( 4.38 , 4.42 )
 };
\addplot[mark=*, very thin] plot coordinates {
 ( 4.56 , 5.1 )
 };
\addplot[mark=*, very thin] plot coordinates {
 ( 4.45 , 5.22 )
 };
\addplot[mark=*, very thin] plot coordinates {
 ( 3.49 , 6.29 )
 };
\addplot[mark=*, very thin] plot coordinates {
 ( 4.23 , 4.34 )
 };
\addplot[mark=*, very thin] plot coordinates {
 ( 4.62 , 5.62 )
 };
\addplot[mark=*, very thin] plot coordinates {
 ( 4.53 , 5.1 )
 };
\addplot[mark=*, very thin] plot coordinates {
 ( 4.45 , 5.22 )
 };
\addplot[mark=*, very thin] plot coordinates {
 ( 4.53 , 5.18 )
 };
\addplot[mark=*, very thin] plot coordinates {
 ( 4.43 , 5.57 )
 };
\addplot[mark=*, very thin] plot coordinates {
 ( 4.38 , 4.62 )
 };
\addplot[mark=*, very thin] plot coordinates {
 ( 4.45 , 5.06 )
 };
\addplot[mark=*, very thin] plot coordinates {
 ( 4.5 , 5.34 )
 };
\addplot[mark=*, very thin] plot coordinates {
 ( 4.45 , 5.34 )
 };
\addplot[mark=*, very thin] plot coordinates {
 ( 4.55 , 5.54 )
 };
\addplot[mark=*, very thin] plot coordinates {
 ( 4.45 , 4.98 )
 };
\addplot[mark=*, very thin] plot coordinates {
 ( 4.42 , 4.5 )
 };

        \end{axis}

\end{tikzpicture}
%\caption{Estimated lines for the stars data in \cite{stars} for the classical models. \label{f:fig3}}
%\end{figure}
}
\end{subfigure}
\caption{Optimal Lines with the classical methods for the stars data set.\label{f:fig3}}
\end{center}
\end{figure}

Observe that the LSS and LAD models were not able to adequately fit the the data while the others (which are somehow similar) show their better performance against  the outliers. Note also that $\GCoD$ reflects this fact, although it is not clear whether LTS(75) (the one with the largest $\GCoD$) is better than the others.

In order to show the behavior of the LTS models and which are the results of their optimal fitting lines, Figure \ref{f:fig3} shows the fitting lines that minimize the 25\%, 50\%, 75\% or 90\% of the residuals and the points that the corresponding optimization problems discard (filled dots in the subfigures) to reach the fitted lines.

\input{figura-ffig4.tex}
\end{ex}
Apart from the classical models described above, the standard vertical distance residuals may be aggregated using a general $\Phi$ function as those introduced  in \eqref{phi}  providing a wide family of new methods to compute the coefficient of the best fitting hyperplanes. Also, linear constrained versions of the above methods may be considered by adding the adequate constraints to the corresponding formulations. Furthermore, many other alternative methods that use vertical distance residuals as MINSADBED  or convex combinations of LSS and LAD methods \cite{arthanary-dodge80} can easily be cast into our modelling framework. The formulations that allow solving those problems are rather similar to those already presented in this section and therefore are left for the interested reader.

\section{Fitting Hyperplanes with block-norm residuals}
\label{sec:3}

In this section we present  models to compute the parameters of fitting hyperplane when the distance point-to-hyperplane is assumed to be a block-norm distance between the point and the closest  point in the hyperplane; and the aggregation criterion is considered in the general form given by \ref{eq1}. Recall that a block norm is a norm such that its unit ball is a polytope symmetric with respect to the origin and with non empty interior.  Block norms, also referred to as polyhedral norms, play an important role in the
measurement of distances in many areas of Operations Research and Applied Mathematics as for instance in Location analysis or Logistics. They are often used to model real world situations
(like measuring highway distances) more accurately than the standard
Euclidean norm. In addition, they can also be used to approximate arbitrary norms since
the set of block norms is dense in the set of all norms \cite{ward-wendell85}.

We denote by $\|\cdot\|_B$ the norm in $\R^d$ whose unit ball is given by a symmetric with respect to the origin, with non empty interior polytope $B$, i.e. $B = \{x \in \R^d: \|x\|_B \leq 1\}$. Let $\ext(B) = \{b_g: g=1, \ldots, G\}$ be the set of extreme points of $B$ and $B^0$ the polar set of $B$ which is defined as:
$$
B^0 = \{v \in \R^d: v^t b_g  \leq 1, g=1, \ldots, G\}
$$
and $\ext(B^0)= \{b^0_1, \ldots, b^0_{G^0}\}$.

The following result characterizes the expression of a block-norm distance in terms of the extreme points of the polar set of the polytope $B$.
\begin{lem}[Ward and Wendell \cite{ward-wendell80,ward-wendell85}]
\label{lem:ward}
Let $B$ be a polytope in $\R^d$ and $x \in \R^d$, then:
$$
\|x\|_B = \max \{|x^t b_g^0|: g=1, \ldots, G^0\}.
$$
\end{lem}
Special cases of block norms are the Manhattan ($\ell_1$) and the Chebyshev ($\ell_\infty$) norms for adequate choices of the extreme points of the unit balls. For instance in $\R^2$, such distances are characterized by the following set of extreme points of their unit balls, $\{ \pm(1, 0), \pm (0,1)\}$ and $\{\pm (1,1), \pm (1,-1)\}$, respectively. Any block norm $\|\cdot\|_B$ in $\R^d$ induces a distance between vectors $x, y \in \R^d$  given by $\d_B(x,y) = \|x-y\|_B$.

Given a set of points $\{x_1, \ldots, x_n\} \subseteq \R^d$ and a  polyhedral unit ball $B$, our goal is to obtain the hyperplane $\H(\bbeta)=\{y \in \R^d: (1,y^t)\bbeta = 0\}$ such that the overall distance $\d_B(\cdot, \cdot)$ from the sample to $\H(\bbeta)$ is minimized according to the globalizing criterion $\Phi$ (for $1 \leq p = \frac{r}{s} \in \Q$). That is:
\begin{equation}
\label{eq:B}\tag{${\rm RM}_B$}
\min_{\bbeta \in \R^{d+1}} \; \dsum_{i=1}^n \lambda_i \bepsilon_{(i)}^{p}
\end{equation}
where for any $x \in \R^d$, $\bepsilon_x=\d_B(x, \H(\bbeta))= \min_{z \in \H(\bbeta)} \d_B(x,z)$, is the ``$\|\cdot\|_B$-projection'' of $x$ onto the hyperplane $\H(\bbeta)$, and $\bepsilon_{(i)}$ denotes the element in $\{\bepsilon_{x_1}, \ldots, \bepsilon_{x_n}\}$ which is sorted in the $i$-th position (in nondecreasing order).

We recall that according to equation \eqref{eq:normdist} in Lemma \ref{le:response}, for any polytope $B$ symmetric with respect to the origin  and with  non empty interior, and $\H(\bbeta)=\{ y^t \in \R^{d}: (1,y^t) \bbeta = 0\}$ then
$\d_{B}(x_{-0}, \H(\bbeta)) = \dfrac{| \bbeta^t x|}{\|\bbeta_{-0}\|_{B^0}},$ where $B^0$ is the polar set of $B$ and $x^t=(1, x_1, \ldots, x_d) \in \R^{d+1}$ is a given point.

\begin{lem}
\label{lem:norma1}
Let $\bbeta^* \in \R^{d+1}$ be an optimal solution of \ref{eq:B} with $\bbeta^*_{-0}  \neq 0$. Then, $\bbeta^\prime = \dfrac{\bbeta^*}{\|\bbeta_{-0}\|_{B^0}}$ is also an optimal solution of \ref{eq:B} with $\|\bbeta^\prime_{-0}\|_{B^0}=1$. Thus, there is an optimal solution of \ref{eq:B}, $\bbeta$, that verifies $\d_B(x_{-0}, \H(\bbeta)) = |\bbeta^t x|$ for any $x^t = (1, x_1, \ldots, x_d) \in \R^{d+1}$.
\end{lem}

From the above lemma, we have
 \begin{theo}
Let $\{x_1, \ldots, x_n\} \subset \R^{d+1}$ be a set of points and let $B \subset \R^d$ be a polytope with $\ext(B)=\{b_1, \ldots, b_G\}$. Then, \ref{eq:B} is equivalent to the following disjunctive programming problem
{\small
\begin{align}
\rho^*(B) :=& \min  \dsum_{j=1}^n \lambda_j \theta_j\label{p:blocknorm}\tag{${\rm LRP}_{\Phi, B}$}\\
\mbox{ s.t. } \mbox{ } \quad & \eqref{ctr:1}-\eqref{ctr:5} \nonumber\\
&\bepsilon_i \geq \bbeta^t x_i, \forall i=1, \ldots, n,\label{B:c1}\\
& \bepsilon_i \geq -\bbeta^t x_i, \forall i=1, \ldots, n,\label{B:c2}\\
&   \bbeta_{-0}^t b_g \leq 1,\; \forall g=1,\ldots,G,\label{B:c4}\\
&   \bigvee_{g=1}^G \bbeta_{-0}^t b_{g}= 1,\label{B:disjunctive}\\
&  \bbeta \in \R^{d+1}, z, \theta, e \in \R^n.\nonumber\\
&  w_{ij} \in \{0,1\}, \forall i, j=1, \ldots, n.\nonumber
\end{align}
}
\end{theo}
\proof
Let us denote by $\bepsilon_i = \d_B(x_i, \H(\bbeta))$. By Lemma \ref{le:response}, $\bepsilon_i = \dfrac{|\bbeta^t x_i|}{\|\bbeta_{-0}\|_{B^0}}$. Furthermore, by Lemma \ref{lem:norma1}, we can assume that $\|\bbeta_{-0}\|_{B^0} = 1$, hence $\bepsilon_i = |\bbeta^t x_i|$ (constraints \eqref{B:c1} and \eqref{B:c2}). By Lemma \ref{lem:ward}, $\|\bbeta_{-0}\|_{B^0} =  \max \{|\dsum_{i=1}^d \beta_i b_{gi}|: g=1, \ldots, G\}$ since $(B^0)^0 = B$. Hence, there exists $g_0 \in \{1, \ldots, G\}$ such that $\|\bbeta_{-0}\|_{B^0} = 1$ (disjunctive constraint \eqref{B:disjunctive}) and thus $  \dsum_{k=1}^d \beta_k b_{gk} \le \dsum_{k=1}^d \beta_k b_{g_0k}=1 $ (constraint \eqref{B:c4}). (Note that absolute values do not need to be taken explicitly into account since if $b_g \in \ext(B)$, then $-b_g \in \ext(B)$.)
\endproof

The above problem can be equivalently written as an  unique mixed integer second order cone programming problem once constraints (\ref{ctr:2}) are transformed using the result in Lemma \ref{le:n2} and binary variables are added to decide which $g_0$ is chosen to verify constraint \eqref{B:c4}. By the same token, this problem can be also equivalently rewritten as  $G$ different SOCP programming problems (each of them fixed to verify one of the disjunctive constraints). Furthermore, MINLP disjunctive programming techniques (e.g. \cite{balas79}, \cite{LeeGrossmann}) may be used to solve the corresponding problem. The following result states a MINLP formulation for \ref{eq:B}:

\begin{cor}
Let $\{x_1, \ldots, x_n\} \subset \R^{d+1}$ be a set of points and let $B \subset \R^d$ be a polytope with $\ext(B)=\{b_1, \ldots, b_G\}$. Then,  \ref{p:blocknorm}  is equivalent to the following problem:
{\small
\begin{align}
\rho^*(B) :=& \min  \dsum_{j=1}^n  \lambda_j  \theta_{j} \label{p:blocknorm}\tag{${\rm LRP}_{\Phi, B}$}\\
\mbox{ s.t. } \mbox{ } \quad & \eqref{ctr:1}-\eqref{ctr:5} \nonumber\\
&\bepsilon_{i} \geq \bbeta_h^t x_i, \forall i=1, \ldots, n,  h=1, \ldots, G,\label{BD:c1}\\
& \bepsilon_{i} \geq -\bbeta_h^t x_i, \forall i=1, \ldots, n,  h=1, \ldots, G,\label{BD:c2}\\
&   \bbeta_{-0h}^t b_g \leq 1,\; \forall g=1,\ldots,G,   h=1, \ldots, G,\label{BD:c4}\\
&   \bbeta_{-0h}^t b_{h} = \xi_h,   h=1, \ldots, G,\label{BD:c3}\\
& \dsum_{h=1}^G \xi_h = 1, \label{BD:c4}\\
&  z, \theta, \bepsilon \in \R^n,\nonumber\\
& \bbeta_h \in \R^{d+1},, \xi_h \in \{0, 1\}, \forall  h=1, \ldots, G,\nonumber\\
&  w_{ij} \in \{0,1\}, \forall i, j=1, \ldots, n.\nonumber
\end{align}
}
\end{cor}

Some special cases for the globalizing criterion $\Phi$  allow even simpler formulations reducing considerably the computational complexity of the problems. In particular, when $\lambda_i=1$ for all $i=1, \ldots, n$, the integer variables  representing ordering ($w_{ij}$)  can be removed from the above formulation.

The following result will allow us to consider polyhedral norms which are \textit{dilations} of other polyhedral norms, i.e., polyhedral norms $\|\cdot\|_{\mu B}$ for some bounded polyhedron $B$ and $\mu >0$ ($\mu B = \{\mu \ z: z \in B\}$).

\begin{cor}
\label{cor:dilation}
Let $\overline{B}$ be a polytope and $\mu >0$. Then, if $\bbeta^*$ is an optimal solution for \ref{p:blocknorm} for $B=\overline{B}$, $\widehat{\bbeta} = \frac{1}{\mu} \bbeta^*$ is an optimal solution for \ref{p:blocknorm} when $B=\mu \overline{B}$. Moreover, $\rho^*(\mu \overline{B})= \frac{1}{\mu^p} \rho^*(\overline{B})  $.
\end{cor}

\proof
It is sufficient to observe that for any $\bbeta \in \R^{d+1}$:
 \begin{align*}
  \|(\beta_1, \ldots, \beta_d)\|_{\mu \overline{B}^0} &= \max \{|\mu b_g^t \bbeta^t|: g=1, \ldots G\}\\
   &= \mu \max \{|b_g^t \bbeta^t|: g=1, \ldots G\} =  \mu\|(\beta_1, \ldots, \beta_d)\|_{\overline{B}^0}.
    \end{align*}

Since $\Phi_{\mu \overline{B}} (\bepsilon_1, \ldots, \bepsilon_n) = \frac{1}{\mu^p} \Phi_{\overline{B}}(\bepsilon_1, \ldots, \bepsilon_n)$, we get the relation between the optimal values. Let $\bbeta^*$ be an optimal solution of \ref{p:blocknorm}. Then, $\dfrac{1}{\mu} \bbeta^*$ is clearly a feasible solution to \ref{p:blocknorm} when $B=\mu \overline{B}$ since $\|(\dfrac{1}{\mu} \beta_1^*, \ldots, \dfrac{1}{\mu} \beta_d^*)\|_{\mu \overline{B}^0} = \|(\beta_1^*, \ldots, \beta_d^*)\|_{\overline{B}^0} = 1$.
\endproof

For the sake of computing $\GCoD$, for solutions to problems with block-norm residuals, note that the one dimensional problem \ref{eq2} does depend on $\Phi$ and also on the residuals through $\kappa_\varepsilon$.  Let us denote by $\kappa_B$ the constant $\kappa_\varepsilon$ when the residuals $\varepsilon_x$ are defined as the block-norm projection with unit ball given by the polytope $B$.

\begin{cor}
Let $B \subset \R^d$ be a polytope. The Goodness of fitting index, $\GCoD$,  when the residuals are defined as the block-norm distance with unit ball $B$, can be computed as:
$$
\GCoD_{\Phi, \varepsilon} = 1- \dfrac{\Phi^* }{\dsum_{i=1}^n |x_{id} - ((\lambda,p)-{\rm omp}(x_{\cdot d}))|^p} \cdot {\dmax_{g=1,\ldots, G} |b_{gd}|},
$$
where $(\lambda,p)$-omp$(x_{\cdot d})$ is the solution to the problem \ref{eq2} with residuals measured with the polyhedral norm with unit ball $B$.
\end{cor}
\proof
By Lemma \ref{le:ompres} the goodness of fitting index $\GCoD_{\Phi,\bepsilon}$ can be computed as:
\begin{equation} \label{eq:det-poly}
\GCoD_{\Phi, \varepsilon} = 1- \dfrac{\Phi^* }{\min_{\beta_0\in \R} \Phi(\kappa_B |x_{1d} - \beta_0|,\ldots, \kappa_B |x_{nd} - \beta_0|)},
\end{equation}
where $\kappa_B=\frac{1}{\dmax_{z\in B} z_d}$.

Observe that since $B$ is a polytope then the above maximum is attained in an extreme point of $B$, namely $b_1,\ldots,b_G$; and thus $\kappa_B=\frac{1}{\dmax_{g=1,\ldots,G} b_{gd}}$.

Next, the problem \ref{eq2} in this case can be expressed as:
$$ \kappa_B \min _{\beta_0\in \R} \sum_{i=1}^n \lambda_i |x_{\cdot d}-\beta_0|_{(i)}^p.$$
Recall that this is a $(\lambda,p)$ \textit{Ordered median problem} and that its optimal solution, a $(\lambda,p)$-omp,  can be easily obtained by the result in Theorem \ref{teor:ompset}.
Replacing the optimal solution to this problem in (\ref{eq:det-poly}) it results in:
$$
 \GCoD_{\Phi, \varepsilon} = 1- \dfrac{\Phi^* }{\dsum_{i=1}^n |x_{id} -( (\lambda,p)-{\rm omp}(x_{\cdot d}))|^p} \cdot {\dmax_{g=1,\ldots, G} |b_{gd}|}.
$$
\endproof

Note that for $\lambda=(1,\ldots,1)$ the $(\lambda,1)$-omp is the standard median point and thus the expression $\dsum_{i=1}^n |x_{id} - {\rm median}(x_{\cdot d})|$ is what it is usually called the \textit{mean absolute deviation with respect to the median}. It is a well-known criterion to find robust optimal hyperplanes of the mean value and a direct measure of the scale of a random variable about its median with  many applications in different fields (see \cite{pham}).

We illustrate the behavior of the block-norm residuals fitting hyperplanes with the same data set used in the Section \ref{sec:2}.

\begin{ex}
\label{ex:stars2}
We consider again the stars data used in Example \ref{ex:stars}. In this case, we run our implementation in \r{} for $\ell_1$-norm, $\ell_\infty$-norm  and hexagonal norm (as the one used in \cite{NP05} with $\ext(B) = \{\pm(2,0), \pm(2,2),\pm (-1,2)\}$) residuals. We use three different criteria: overall SUM ($\lambda=(1,\ldots, 1)$ and $p=1$), MAXimum ($\lambda=(1,0,\ldots, 0)$ and $p=1$),  $K$-centrum ($\lambda=(\overbrace{0, \ldots, 0}^{K},  \overbrace{1, \ldots, 1}^{n-K})$) for $K=\lfloor 0.75n\rfloor$ (the model will minimize the sum of the $25\%$ greatest residuals) and anti-$K$-centrum ($\lambda=(\overbrace{1, \ldots, 1}^{K},  \overbrace{0, \ldots, 0}^{n-K})$) for $K=\lfloor 0.5n\rfloor$ (the model will minimize the sum of the $50\%$ smallest residuals).  The results for all the combinations and the graph for the $K$-centrum lines are shown in Figure \ref{f:fig-star}.

%\begin{table}[h]
%\centering
\begin{figure}[t!]
\centering
\begin{subfigure}{0.6\linewidth}
\scalebox{0.9}{\begin{tabular}{ccc}
Method $(\Phi, \bepsilon)$ & Optimal Line & $\GCoD_{\Phi, \bepsilon}$\\\hline\hline
(SUM, $\ell_1$) &  $y = 7 x  -25.81$ & $0.6505853$\\
(SUM, $\ell_\infty$) & $y = 5.25 x + -18.1425$ & $0.7009688$ \\
(SUM, Hex) &  $y = 7 x  -25.81$ & $0.6505853$\\\hline
(MAX, $\ell_1$) &  $y= -3.230769 x+ 18.77577$ & $0.5336373$\\
(MAX, $\ell_\infty$) & $y= -3.230769 x+ 18.77577$ & $0.6438685$ \\
(MAX, Hex)& $y= -3.230769 x+ 18.77577$ & $0.6438685$ \\\hline
(kC, $\ell_1$) &  $y= -4.307692 x + 23.03346$ & $0.4628481$\\
(kC, $\ell_\infty$) & $y= -2.493333 x +15.67113$ & $0.5921635$ \\
(kC, Hex) &  $y= 7.642857x+  -28.67929$  & $0.8317972 $\\\hline
(AkC, $\ell_1$) &  $y= 5.6 x -19.804$ & $0.8443055$\\
(AkC, $\ell_\infty$) & $y= 4.869565 x -16.41565$ & $0.8426523$ \\
(AkC, Hex) &  $y=5.473684 x -19.28316$  & $0.6431602$\\\hline
\end{tabular}}
\end{subfigure}
\hspace*{-0.7cm}
\begin{subfigure}{0.4\linewidth}
\scalebox{0.9}{%\begin{figure}[h]
%\centering
\begin{tikzpicture}[scale=0.55]
  \begin{axis}[
  xmin=3.4,   xmax=5.0,
	ymin=4.0,   ymax=6.5,
xtick = {3.5, 3.7, 3.9, 4.1, 4.3, 4.5, 4.7, 4.9},
ytick = {4.0, 4.5, 5.0, 5.5, 6.0, 6.5},
    %anchor=origin,  % Align the origins
    x=7cm, y=4cm,   % Set the same unit vectors
    ]

         \addplot[mark=none, red, thick] plot coordinates {
        (3.4, -4.307692*3.4 + 23.03346)
                (4.8, -4.307692*4.8 + 23.03346)
    };

     \addplot[mark=none, blue, thick] plot coordinates {
        (3.4,-2.493333*3.4 +15.67113)
                (4.8,-2.493333*4.8 +15.67113)
    };

      \addplot[mark=none, black, dotted, ,thick] plot coordinates {
        (3.4,7.642857*3.4+  -28.67929)
                (4.8, 7.642857*4.8+  -28.67929)
    };

       \addplot[mark=none, green, thick] plot coordinates {
        (3.4,-0.413303*3.4 + 6.793467)
                (4.8,-0.413303*4.8+ 6.793467)
    };

    \legend{$\ell_1$, $\ell_\infty$, Hex, LSS}

\addplot[mark=*, very thin] plot coordinates {
 ( 4.37 , 5.23 )
 };
\addplot[mark=*, very thin] plot coordinates {
 ( 4.56 , 5.74 )
 };
\addplot[mark=*, very thin] plot coordinates {
 ( 4.26 , 4.93 )
 };
\addplot[mark=*, very thin] plot coordinates {
 ( 4.56 , 5.74 )
 };
\addplot[mark=*, very thin] plot coordinates {
 ( 4.3 , 5.19 )
 };
\addplot[mark=*, very thin] plot coordinates {
 ( 4.46 , 5.46 )
 };
\addplot[mark=*, very thin] plot coordinates {
 ( 3.84 , 4.65 )
 };
\addplot[mark=*, very thin] plot coordinates {
 ( 4.57 , 5.27 )
 };
\addplot[mark=*, very thin] plot coordinates {
 ( 4.26 , 5.57 )
 };
\addplot[mark=*, very thin] plot coordinates {
 ( 4.37 , 5.12 )
 };
\addplot[mark=*, very thin] plot coordinates {
 ( 3.49 , 5.73 )
 };
\addplot[mark=*, very thin] plot coordinates {
 ( 4.43 , 5.45 )
 };
\addplot[mark=*, very thin] plot coordinates {
 ( 4.48 , 5.42 )
 };
\addplot[mark=*, very thin] plot coordinates {
 ( 4.01 , 4.05 )
 };
\addplot[mark=*, very thin] plot coordinates {
 ( 4.29 , 4.26 )
 };
\addplot[mark=*, very thin] plot coordinates {
 ( 4.42 , 4.58 )
 };
\addplot[mark=*, very thin] plot coordinates {
 ( 4.23 , 3.94 )
 };
\addplot[mark=*, very thin] plot coordinates {
 ( 4.42 , 4.18 )
 };
\addplot[mark=*, very thin] plot coordinates {
 ( 4.23 , 4.18 )
 };
\addplot[mark=*, very thin] plot coordinates {
 ( 3.49 , 5.89 )
 };
\addplot[mark=*, very thin] plot coordinates {
 ( 4.29 , 4.38 )
 };
\addplot[mark=*, very thin] plot coordinates {
 ( 4.29 , 4.22 )
 };
\addplot[mark=*, very thin] plot coordinates {
 ( 4.42 , 4.42 )
 };
\addplot[mark=*, very thin] plot coordinates {
 ( 4.49 , 4.85 )
 };
\addplot[mark=*, very thin] plot coordinates {
 ( 4.38 , 5.02 )
 };
\addplot[mark=*, very thin] plot coordinates {
 ( 4.42 , 4.66 )
 };
\addplot[mark=*, very thin] plot coordinates {
 ( 4.29 , 4.66 )
 };
\addplot[mark=*, very thin] plot coordinates {
 ( 4.38 , 4.9 )
 };
\addplot[mark=*, very thin] plot coordinates {
 ( 4.22 , 4.39 )
 };
\addplot[mark=*, very thin] plot coordinates {
 ( 3.48 , 6.05 )
 };
\addplot[mark=*, very thin] plot coordinates {
 ( 4.38 , 4.42 )
 };
\addplot[mark=*, very thin] plot coordinates {
 ( 4.56 , 5.1 )
 };
\addplot[mark=*, very thin] plot coordinates {
 ( 4.45 , 5.22 )
 };
\addplot[mark=*, very thin] plot coordinates {
 ( 3.49 , 6.29 )
 };
\addplot[mark=*, very thin] plot coordinates {
 ( 4.23 , 4.34 )
 };
\addplot[mark=*, very thin] plot coordinates {
 ( 4.62 , 5.62 )
 };
\addplot[mark=*, very thin] plot coordinates {
 ( 4.53 , 5.1 )
 };
\addplot[mark=*, very thin] plot coordinates {
 ( 4.45 , 5.22 )
 };
\addplot[mark=*, very thin] plot coordinates {
 ( 4.53 , 5.18 )
 };
\addplot[mark=*, very thin] plot coordinates {
 ( 4.43 , 5.57 )
 };
\addplot[mark=*, very thin] plot coordinates {
 ( 4.38 , 4.62 )
 };
\addplot[mark=*, very thin] plot coordinates {
 ( 4.45 , 5.06 )
 };
\addplot[mark=*, very thin] plot coordinates {
 ( 4.5 , 5.34 )
 };
\addplot[mark=*, very thin] plot coordinates {
 ( 4.45 , 5.34 )
 };
\addplot[mark=*, very thin] plot coordinates {
 ( 4.55 , 5.54 )
 };
\addplot[mark=*, very thin] plot coordinates {
 ( 4.45 , 4.98 )
 };
\addplot[mark=*, very thin] plot coordinates {
 ( 4.42 , 4.5 )
 };

        \end{axis}

\end{tikzpicture}
%\caption{Estimated lines for the Hertzsprung--Russell stars data in Example \ref{ex:stars2} for the $K$-centrum criteria compared with the LSS line. \label{f:fig-star}}
%\end{figure}
}
\end{subfigure}
\caption{Optimal lines obtained with block-norm residuals for the stars data set.\label{f:fig-star}}
\end{figure}
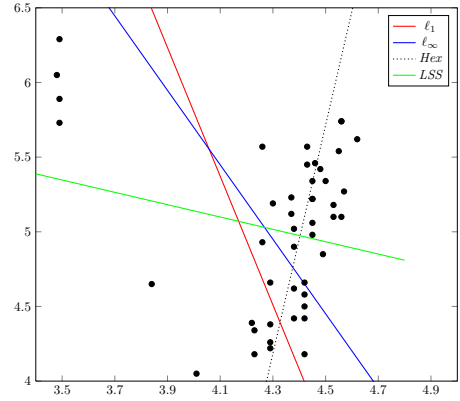

Note that different situations may happen when running the different models: in the case of the SUM criterion the models for $\ell_1$ and hexagonal residuals coincide; in the MAX criterion the three optimal lines are the same, and for the $K$-centrum and anti-$K$-centrum the three models are different. Furthermore, even in the case when the models coincide, one may have different goodness of fitting indices due to the different way of measuring distances (see the $\ell_1$ and hexagonal residuals for the MAX criterion).

From the above, we observed that the $\GCoD$ are not comparable when different residuals are used in the models since the value given to the residuals (both with respect to the best model and with respect to the simplified model with only intercept) is different. Thus, the generalized coefficient allows us to compare the goodness of fitting between models provided that the distance (to measure the residuals) and the aggregation criterion are fixed.

\end{ex}

\section{Fitting Hyperplanes with $\ell_\tau$ distances}
\label{sec:4}
In this section we present the mathematical programming formulations for computing the optimal hyperplanes when the residuals are defined as $\ell_\tau$ distances between demand points and the linear body. Recall that the $\ell_\tau$-norm in $\R^d$,  with $\tau\geq 1$, is defined as:
$$
\|z\|_\tau = \left\{ \begin{array}{ll} \left( \dsum_{k=1}^d |z_k|^\tau \right)^{\frac{1}{\tau}} & \mbox{$if \tau<\infty$},\\
\dmax_{k=1, \ldots, d} \{|z_k|\} & \mbox{if $\tau=\infty$}\end{array}\right.
$$
for any $z= (z_1, \ldots, z_d)^t \in \R^d$. From this norm we denote by $\d_{\ell_\tau}(z,y) = \|z-y\|_\tau$ the $\ell_\tau$-distance between the points $z, y \in \R^d$. The well-known Euclidean distance that measures the straight line distance between points in $\R^d$ is the $\ell_2$-norm in this family.  Note that the extreme cases of $\ell_1$ and $\ell_\infty$  represent both  block   and $\ell_\tau$-norms, since their unit balls are polytopes but also fit within the family of $\ell_\tau$-norms.

We recall that according to equation \eqref{eq:normdist} in Lemma \ref{le:response}, for any $\tau = \frac{r}{s} \in \Q$ with $r\geq s \in \Z_+$, $\gcd(r,s)=1$ and $\H(\bbeta)=\{ y^t \in \R^{d}: (1,y^t) \bbeta = 0\}$ then
$\d_{\tau}(z, \H(\bbeta)) = \dfrac{| \bbeta^t z|}{\|\bbeta_{-0}\|_\nu}
$
where $\nu$ is such that $\dfrac{1}{\tau} + \dfrac{1}{\nu} = 1$ (for $\tau=1$, $\nu=\infty$ while for $\tau=\infty$, $\nu=1$).% = \left\{\begin{array}{rl}  \dfrac{\tau}{\tau-1} = \frac{r}{r-s} & \mbox{if $r\neq s$,}\\
% \infty & \mbox{if $r=s$ ($\tau=1$),}\\
 %1 & \mbox{if $\tau=\infty$}\end{array}\right.$.

In this section we will assume that the residuals are defined as  the shortest distance from the points to the fitted hyperplane, namely to their projections, under a given $\ell_\tau$ norm. In other words, for a given point $\hat x=(1,\hat x_1,\dots, \hat x_d)^t$ the residual is:
$\varepsilon_{\hat x} (\bbeta) = \d_{\tau}(\hat x_{-0}, \H(\bbeta)).$

As in previous sections, for a given set of points $\{x_1, \ldots, x_n\} \subseteq \R^{d+1}$, the computation of the parameters $\bbeta \in \R^{d+1}$  assuming that the globalizing criterion is $\Phi$ and the residuals are measured with $\ell_\tau$-distance can be obtained by solving an adequate optimization problem.

\begin{theo}
Let $\{x_1, \ldots, x_n\} \subset \R^{d+1}$ be a set of points, $\lambda \in \R^n$, $\tau = \dfrac{r}{s}\in \Q$ with $r > s \in \N$ and $\gcd(r,s)=1$, and $\|\cdot\|_\tau$ a $\ell_\tau$-norm in $\R^d$. The Problem \ref{eq1} is equivalent to the following mathematical programming problem:

\begin{align}
\min \dsum_{j=1}^n \lambda_j \theta_j\label{eq:lpres}\tag{${\rm LPR}_{\Phi, \ell_\tau}$}\\
\mbox{ s.t. } \mbox{ } \quad & \eqref{ctr:1}-\eqref{ctr:5}, \eqref{B:c1}-\eqref{B:c2}, \nonumber\\
& \|\bbeta_{-0}\|_\nu = 1, \label{normctr}\\
& w_{ij} \in \{0,1\}, & \forall i=1, \ldots, n,\nonumber\\
& z, \theta \in \R^n_+, \bbeta \in \R^{d+1}.\nonumber
\end{align}

%with $\nu=\frac{\tau}{\tau-1}$.

\end{theo}

Note that the above problem is nonconvex for $1<\tau<\infty$ because of the binary variables and constraint \eqref{normctr}. Approximation schemes are available in different free and commercial solvers, although no guarantee of optimality is provided (e.g., NLOPT, MATLAB, Minotaur, ...). In what follows we describe an approximation approach based on some linear approximations of the problem.

Let $P$ be a polyhedron such that $P \subset \mathcal{B}= \{z \in \R^d: \|z\|_\nu \leq 1\}$, and denote by $r_P= \sup_{\|z\|_P = 1} \|z\|_\nu$ (note that by construction $r_P \leq 1$). Observe that $r_P$ is the radius of the smallest $\ell_\nu$-ball containing $P$. In addition, let $Q$ be a polyhedron such that $\mathcal{B} \subset Q$, and denote by $R_Q= \inf_{\|z\|_Q = 1} \|z\|_\nu$ (note that by construction $R_Q \geq 1$). In this case $R_Q$ is the radius of the largest $\ell_\nu$-ball contained in $Q$.

\begin{theo}
\label{theo:P}
Let $\lambda_1, \ldots, \lambda_n \geq 0$  and the globalizing function $\Phi(\bepsilon_1, \ldots, \bepsilon_n) = \dsum_{i=1}^n \lambda_i \bepsilon_{(i)}^{\delta}$ then:
\begin{eqnarray}
\Phi_{P^*} \leq \Phi_{\ell_\tau} \leq \dfrac{1}{r_P^{ \delta}}\Phi_{P^*} \label{eq:rP}\\
\frac{1}{R_Q^{ \delta}} \Phi_{Q^*} \leq \Phi_{\ell_\tau} \leq \Phi_{Q^*} \label{eq:RQ}
\end{eqnarray}
\end{theo}

\proof

By the relations between the norms, it is clear that $\|z\|_P \geq \|z\|_\nu \geq r_P \|z\|_P$. Let $\H(\bbeta)=\{z \in \R^d: (1,z^t)\bbeta =0\}$. Then, for any  $x \in \R^d$, the above relationships imply the following inequalities relating the distances with respect to $\|\cdot\|_{P^*}$-residuals and $\|\cdot\|_{\tau}$-residuals:

$$
\d_{P^*}(x_{-0}, \H(\bbeta)) = \dfrac{|\bbeta^t x|}{\|\bbeta_{-0}\|_P} \leq \dfrac{|\bbeta^t x|}{\|\bbeta_{-0}\|_\nu} \leq \, \d_{\tau}(x_{-0}, \H(\bbeta))
$$

and

$$
\d_{\tau}(x_{-0}, \H(\bbeta)) = \dfrac{|{\bbeta}^t x|}{\|\bbeta_{-0}\|_\nu} \leq \dfrac{|{\bbeta}^t x|}{r_P\|\bbeta_{-0}\|_{P}} \leq \frac{1}{r_P} \d_{P^*}(x_{-0}, \H(\bbeta))
$$

Let us consider the globalizing  criterion $\Phi(\bepsilon_1, \ldots, \bepsilon_n) = \dsum_{i=1}^n \lambda_i \bepsilon_{(i)}^{ \delta}$. Then, the evaluation of $\Phi$ with respect to the residuals computed with the polyhedral norm with unit ball  $P^*$ and the $\ell_\tau$-norm, namely $\varepsilon_{i, P^*}= \d_{P^*}(x_{i,-0}, \H(\bbeta))$ and $\varepsilon_{i,\ell_\tau}= \d_{\tau}(x_{i,-0}, \H(\bbeta))$ for all $i=1, \ldots, n$, satisfies:
$$
\Phi(\bepsilon_{P^*}) \leq \Phi(\bepsilon_{\ell_\tau}) \leq \dfrac{1}{r_P^{ \delta}} \Phi(\bepsilon_{P^*}).
$$
This equation proves (\ref{eq:rP}).

Next, it is clear that $\|z\|_Q \leq \|z\|_\nu \leq R_Q \|z\|_Q$. Now, using an argument similar to the one above we conclude  that
\begin{align*}
\d_{Q^*}(x_{-0}, \H(\bbeta)) &= \dfrac{|\bbeta^t x|}{\|\bbeta_{-0}\|_Q} \geq \dfrac{|\bbeta^t x|}{\|\bbeta_{-0}\|_\nu} \geq \d_{\tau}(x_{-0}, \H(\bbeta))\\
 &= \dfrac{|{\bbeta}^t x|}{\|\bbeta_{-0}\|_\nu} \geq \dfrac{|{\bbeta}^t x|}{{R_Q}\|\bbeta_{-0}\|_\nu} \geq \frac{1}{R_Q} \d_{Q^*}(x_{-0}, \H(\bbeta)).
\end{align*}
From these inequalities it clearly follows (\ref{eq:RQ}).
\endproof

Let $P_N$ be  a symmetric with respect to the origin polytope with $N$ vertices, $\{p_1, \ldots, p_N\}$, inscribed in the $\ell_\nu$ hypersphere $\mathcal{B}=\{z \in \R^d: \|z\|_\nu=1\}$ and let $r_{P_N}$ be the radius of the smallest $\ell_\nu$ ball centered at the origin containing  $P_N$. Let $R_{Q_N} = \frac{1}{r_{P_N}}$ and denote by $Q_N$ the $R_{Q_N}$-dilation of $P_N$. By construction $P_N \subset \mathcal{B} \subset Q_N$. Hence, { for the globalizing function $\Phi(\bepsilon_1, \ldots, \bepsilon_n) = \dsum_{i=1}^n \lambda_i \bepsilon_{(i)}^{\delta}$, by the Theorem \ref{theo:P},}  we get that:

$$
\max\{\Phi(\bepsilon_{P_N^*}), \frac{1}{R_{Q_N}^{ \delta}}\Phi(\bepsilon_{Q_N^*}))\} \leq \Phi(\bepsilon_{\ell_\tau}) \leq \min\{ \Phi(\bepsilon_{Q_N^*}), \frac{1}{r_{P_N}^{ \delta}}\Phi(\bepsilon_{P_N^*})\}
$$

Furthermore, by Corollary \ref{cor:dilation}, since $Q_N$ is a dilation of $P_N$, both problems have the same optimal solutions and $\Phi(\bepsilon_{P_N^*}) = r_P^{ \delta} \Phi(\bepsilon_{Q_N^*})$. Hence,
$$
\frac{1}{r_{P_N}^{ \delta}} \Phi(\bepsilon_{P_N^*}) \leq \Phi(\bepsilon_{\ell_\tau}) \leq \Phi(\bepsilon_{Q_N^*}).
$$

It is clear from its definition that $r_{P_N}$ gives the approximation error whenever a $\ell_\nu$-norm is replaced by a polyhedral norm with unit ball $P_N$. This measure can be   explicitly computed  from the set of inequalities that describe the polyhedron.

\begin{lem}
Let $P = \{z \in \R^d: a_ix \le b_i, i=1, \ldots, N\}$ be a polytope, then:
$$
r_{P}=\dmax_{i=1, \ldots, N} \dfrac{b_i}{\|a_i\|_\tau}.
$$
\end{lem}

\proof
First, note that $r_{P}= \sup_{\|z\|_P = 1} \|z\|_\nu = \dmax_{\|z\|_P = 1} \|z\|_\nu$ by the compactness of $P$. Thus, $r_P$ is the $\ell_\nu$-inradius of $P$. Next, by \cite{mangasarian}, the radius of a $\ell_\nu$ ball centered at the origin and reaching the facet $\{x \in \R^d:  a_i^t x \leq b\}$ of $P$ is the $\ell_\nu$ projection of the origin onto that facet, namely $\dfrac{|b_i|}{\|a_i\|_\tau}$.
Hence,  $r_P$ is the maximum of those distances among the $N$ facets defining $P$.
\endproof

\begin{theo}
Let $\{x_1, \ldots, x_n\} \subset \R^{d+1}$ be a set of demand points, $\lambda \in \R_+^n$, $\tau= \dfrac{r}{s}\in \Q$ with $r > s \in \N$,  $\gcd(r,s)=1$ { and the globalizing function $\Phi(\bepsilon_1, \ldots, \bepsilon_n) = \dsum_{i=1}^n \lambda_i \bepsilon_{(i)}^p$.} The following problem provides a   lower bound for Problem \ref{eq:lpres}.

\begin{align}
\rho^* :=&  \min  \dsum_{j=1}^n \lambda_j \theta_j\label{eq:lpin}\tag{Inner-${\rm \ell_\tau}$}\\
 \mbox{ s.t. } &   (\ref{ctr:1})-(\ref{ctr:5}) & \nonumber\\
 & \bepsilon_i \geq |\bbeta^t x_i|, &\forall i=1, \ldots, n,\label{nlc2}\\
&  \|\bbeta_{-0}\|_{P_N} = 1, \label{normctrP}\\
& w_{ij} \in \{0,1\}, & \forall i=1, \ldots, n,\nonumber\\
&  z, \theta \in \R^n_+, \bbeta \in \R^{d+1}\nonumber
\end{align}

Furthermore, $\rho^* \leq \Phi^*_{\ell_\tau} \leq \frac{1}{r_P^p} \rho^*$.

\end{theo}

\begin{cor}
For any data set $\{x_1,\dots,x_n\}\subset \R^{d+1}$ and any $\ell_\tau$-norm with $1<\tau<+\infty$ there exists a polyhedral norm $\| \cdot \|_B$ whose unit ball $B$ has at most $2n$ extreme points and such  that the optimal values of Problem \ref{eq:lpres} and  \ref{p:blocknorm} coincide.
\end{cor}

In \cite{love-morris72} the authors propose a measure of the goodness of approximating a given norm by another norm. This measure was defined in order to quantify the approximation errors when modeling road distances between cities. We redefine this measure to evaluate the approximation errors when approximating $\ell_\tau$ norms via polyhedral norms:
$$
{\rm SD} = \mathop{\dsum_{i=1}^n}_{\d_\tau(x_i, \bbeta) >0} \dfrac{(\d_\tau(x_i, \bbeta)-\d_P(x_i, \bbeta))^2}{\d_\tau(x_i, \bbeta)}
$$

\begin{ex}
\label{ex:stars3}
Let us consider again the stars data from Example \ref{ex:stars}. We run now the models using as aggregation criteria the overall sum of the residuals ($\Phi=SUM$) and the residuals are the $\ell_\tau$ projections of the points onto the optimal line, for $\tau \in \{1.5, 2, 3\}$. The obtained estimations for the  aggregation criterion $\Phi=SUM$ and their goodness of fitting ($\GCoD_{\Phi, \bepsilon}$) are shown in Table \ref{t:table-appro2}. The obtained lines are drawn in Figure \ref{f:figlp2}.

\input{figura-ffiglp2.tex}

 Observe that for this data set, getting high accuracy for the $\ell_\tau$-norm residual problems is possible using small number of vertices ($N$) in the approximation by polyhedral norms. As expected, increasing the number of vertices improves the accuracy at the price of increasing the computation times.

We also computed the optimal lines for different aggregation criteria ($\Phi\in \{\mbox{SUM, MAX, kC, AkC}\}$) with $\ell_\tau$ residuals, $\tau \in \{1.5, 2, 3\}$, using the polyhedral approximation approach with $N=480$ vertices. The results are shown in Table \ref{t:tableappro3}.  The reader may  observe from these results  that the approximation error, although tiny, depends both of the chosen residuals and  aggregation criteria.

\setlength\tabcolsep{4pt}
\begin{table}[h]
\scriptsize\centering
\begin{tabular}{c|c||c|c|c}
\multicolumn{2}{c||}{}  &  $\ell_{1.5}$  & $\ell_{2}$  &  $\ell_{3}$ \\\hline\hline
\multirow{3}{*}{SUM}& Line & $y=5.92 x  -21.1016$ & $y= 6.75 x-24.6975$ & $y=7 x  -25.81$\\
&$\GCoD$  & $0.6643$ & $0.6542$ & $0.6509$ \\
&${\rm SD}$  & $3.36\times 10^{-10}$  & $ 1.73\times 10^{-10} $  &  $1.65\times 10^{-9}$\\\hline
\multirow{3}{*}{MAX}& Model & $y= -3.2307 x + 18.7757$& $y= -3.2307 x + 18.7757 $ & $y= -3.2307 x + 18.7757$ \\
&$\GCoD$  &$0.5805$   & $0.5544$ & $0.5381$ \\
&${\rm SD}$  & $4.07\times 10^{-14}$ & $1.90\times 10^{-12}$ & $3.85\times 10^{-13}$ \\\hline
\multirow{3}{*}{kC}& Model & $y= -2.8133x + 16.9367$ & $y= -3.1756x +18.5100$ & $y= -4.3076 x + 23.0334$\\
&$\GCoD$  & $0.5111$ & $0.4790$ & $ 0.4650$\\
&${\rm SD}$  & $3.51\times 10^{-13}$  & $7.53\times 10^{-10}$ & $9.70\times 10^{-10}$\\\hline
\multirow{3}{*}{AkC}& Model & $y= 6.75 x -25.0875$ & $y = 6.5555 x -24.1533$ & $y= 5.175 x - 17.7146$\\
&$\GCoD$  & $0.8092$ & $0.82512$ & $0.8217$\\
&${\rm SD}$  & $7.15\times 10^{-10}$ & $2.10\times 10^{-9}$ & $5.49\times 10^{-10}$\\\hline
\end{tabular}
\caption{Optimal lines for different criteria and $\ell_\tau$ residuals of Example \ref{ex:stars3}.\label{t:tableappro3}}

\end{table}

Finally, we compare our approximation scheme for $\ell_\tau$ residuals, on this data set,  with other available implementations. Orthogonal Distance Regression (ODR) is a particular case of our general framework where $\ell_2$ residuals are chosen and $\Phi$ is the sum of squares aggregation criterion (note that both approaches coincide when the coefficient of the dependent  coordinate is non zero while such an assumption is not imposed in our models). The package \texttt{pracma} in \r{} allows to compute ODR by using an approximated iterative procedure (see \cite{boggsetal90}). The models obtained with both approaches are shown in the following table, were one can observe that, for this data set, our approach to approximate $\ell_\tau$ distances by polyhedral norms (with $N=320$ vertices) has a better performance on the global error measure of the models (although as expected the models obtained by both methods are almost the same):

\begin{center}
\begin{tabular}{c|c|c}
& ODR  & SOS-$\ell_2$  ({\rm SD}=$9.93 \times 10^{-11}$)\\\hline
Model & $y=-7.05736x + 35.42935$ & $y= -7.098062 x + 35.60477$\\
Global Residuals & $3.959383$ & $3.662783$ \\
\end{tabular}
\end{center}

\end{ex}
\setlength\tabcolsep{6pt}

\section{Experiments}
\label{sec:5}

We tested the proposed models for different data sets in order to show the applicability and the differences of some of the methods detailed in the sections above. Our formulations have been coded in Gurobi 6.0 under \r{} and executed in a PC with an Intel Core i7 processor at 2x 2.40 GHz  and 4 GB of RAM.  As far as we know, the battery of experiments that we performed has never been considered in the literature, since we have compared 42 different methods (several combinations of aggregation criteria and residuals measures).

\subsection{Synthetic Experiments}

We consider a set of randomly generated points with different peculiarities in order to test and compare the described methodologies, following similar schemes that those proposed in \cite{BM14}. We generated $n=100$ data points in dimension $d  \in \{2, 4\}$, $\{x_1, \ldots, x_n\} \subseteq \R^{d+1}$ as follows. Each $x_{ik}$ follows an independent and identically distributed Gaussian distribution with mean $0$ and standard deviation $100$. We fix $\bbeta^t = (0, 1, \ldots, 1) \in \R^{d+1}$. The last coordinate, $x_d$, is chosen as the response and we generate it as:

$
x_{id} = - \dsum_{k=1}^{d-1} x_{ik} + u_i, \qquad \forall i=1, \ldots, n,
$

\noindent where $u_i$ is also  generated as a Gaussian distribution with mean $0$ and standard deviation $10$.

Then, $15\%$ of the data are now corrupted by adding an extra Gaussian term (with mean $0$ and standard deviation $500$) to: (1)  all the components except the last one or (2) to the last coordinate.

For each one of the generated data sets, we run the models that results from the combination of the following aggregation criteria and residuals detailed in Table \ref{tablef}.

\begin{table}
\centering
\begin{tabular}{|c|c|c|c|}\cline{1-2}\cline{4-4}
\multicolumn{2}{|c|}{\bf Aggregation criteria} & & {\bf Residuals}\\\cline{1-2}\cline{4-4}
SUM & $\dsum_{i=1}^n \bepsilon_i$ & & V \\\cline{1-2}\cline{4-4}
MAX  & $ \dmax_{i=1, \ldots, n} \bepsilon_i$ && $\ell_1$\\\cline{1-2}\cline{4-4}
MED  & $ {\rm median} (\bepsilon_1, \ldots, \bepsilon_n)$ && $\ell_\infty$ \\\cline{1-2}\cline{4-4}
kC & $\dsum_{i=1}^{\lfloor 0.5 n\rfloor} \bepsilon_{(i)}$ & & $\ell_{\frac{3}{2}}$ \\\cline{1-2}\cline{4-4}
AkC & $\dsum_{i=\lfloor 0.5 n\rfloor + 1}^n \bepsilon_{(i)}$ && $\ell_{2}$\\\cline{1-2}\cline{4-4}
SOS & $\dsum_{i=1}^n \bepsilon_i^2$ && $\ell_{3}$\\\cline{1-2}\cline{4-4}
1.5SUM & $\dsum_{i=1}^n \bepsilon_i^\frac{3}{2}$ & \multicolumn{2}{c}{}\\\cline{1-2}
\end{tabular}
\caption{Combinations of chosen aggregation criteria and residuals.\label{tablef}}
\end{table}

\setlength\tabcolsep{1pt}

\begin{table}
\caption{Results for bidimensional experiments corrupting the $X$ variables.\label{table:experiments2}}
\centering
{\scriptsize\begin{tabular}{|c|c|c|c|c|}\cline{3-5}
       \multicolumn{2}{c|}{}           & V     & $\ell_1$ & $\ell_\infty$ \\\hline\cline{3-5}
\multirow{4}{*}{SUM} &  $\widehat{\bbeta}$  & $(-1.9587, 0.3011, 1)$   &  $(1.9587, -0.3011, -1)$   &  $(0.4240, -0.9403, -1)$   \\
%      &  $\Phi^*$  & $7070.199$  &  $7070.199$  &  $3854.186$  \\
      &  $\GCoD$  & $0.1456$  &  $0.1456$  &  $0.5342$  \\
      &  $\%$  & $8\%$  &  $8\%$  &  $65\%$  \\
      &  $\epsilon_{90}$  & $141.2995$  &  $141.2995$  &  $87.0871$ \\ \hline
\multirow{4}{*}{MAX} &  $\widehat{\bbeta}$  & $(10.9038, 0.1571, 1)$   &  $(10.9038, 0.1571, 1)$   &  $(10.9038, 0.1571, 1)$   \\
%      &  $\Phi^*$  & $209.7352$  &  $209.7352$  &  $181.2465$  \\
      &  $\GCoD$  & $0.1484$  &  $0.1484$  &  $0.2641$  \\
      &  $\%$  & $10\%$  &  $10\%$  &  $10\%$  \\
      &  $\epsilon_{90}$  & $158.9295$  &  $158.9295$  &  $158.9295$ \\ \hline
\multirow{4}{*}{SOS} &  $\widehat{\bbeta}$  & $(-3.1753, 0.1860, 1)$   &  $(3.1753, -0.1860, -1)$   &  $(-1.8549, 0.2858, 1)$   \\
%      &  $\Phi^*$  & $821349.6$  &  $821349.6$  &  $538587.6$  \\
      &  $\GCoD$  & $ 0.2261$  &  $0.2261$  &  $0.4925$  \\
      &  $\%$  & $8\%$  &  $8\%$  &  $9\%$  \\
      &  $\epsilon_{90}$  & $157.7177$  &  $157.7177$  &  $143.1279$  \\ \hline
\multirow{4}{*}{1.5SUM} &  $\widehat{\bbeta}$  & $(-3.5386, 0.2112, 1)$   &  $(3.5397, -0.2112, -1)$   &  $(0.3967, -0.4136, -1)$   \\
%      &  $\Phi^*$  & $73890.31$  &  $73890.65$  &  $49637.73$  \\
      &  $\GCoD$  & $0.1812$  &  $0.1812$  &  $0.4499$  \\
      &  $\%$  & $8\%$  &  $8\%$  &  $8\%$  \\
      &  $\epsilon_{90}$  & $152.361$  &  $152.3626$  &  $127.4389$  \\ \hline
\multirow{4}{*}{kC} &  $\widehat{\bbeta}$  & $(-3.0188, 0.2328, 1)$   &  $(-3.0188, 0.2328, 1)$   &  $(0.3503, 0.9091, 1)$   \\
%      &  $\Phi^*$  & $5695.716$  &  $5695.716$  &  $3716.344$  \\
      &  $\GCoD$  & $0.1226$  &  $0.1226$  &  $0.4275$  \\
      &  $\%$  & $8\%$  &  $8\%$  &  $60\%$  \\
      &  $\epsilon_{90}$  & $150.5599$  &  $150.5599$  &  $85.1974$  \\ \hline
\multirow{4}{*}{AkC} &  $\widehat{\bbeta}$  & $(5.8180, 0.7718, 1)$   &  $(2.2956, 0.7734, 1)$   &  $(2.6795, 0.9874, 1)$   \\
%      %&  $\Phi^*$  & $699.9667$  &  $701.7374$  &  $176.3393$  \\
      &  $\GCoD$  & $0.6735$  &  $0.9040$  &  $0.9758$  \\
      &  $\%$  & $29\%$  &  $34\%$  &  $70\%$  \\
      &  $\epsilon_{90}$  & $77.4723$  &  $74.8420$  &  $92.8187$  \\ \hline
\multirow{4}{*}{MED} &  $\widehat{\bbeta}$  & $(6.1846, 0.7795, 1)$   &  $(6.1842, 0.7795, 1)$   &  $(1.3314, 0.9890, 1)$   \\
      %&  $\Phi^*$  & $19.15404$  &  $19.1536$  &  $3.777488$  \\
      &  $\GCoD$  & $0.7021$  &  $0.8690$  &  $0.9741$  \\
      &  $\%$  & $31\%$  &  $31\%$  &  $70\%$  \\
      &  $\epsilon_{90}$  & $78.4775$  &  $78.4772$  &  $91.9773$ \\ \hline
 \multicolumn{5}{c}{}\\\cline{3-5}
      \multicolumn{2}{c|}{}      & $\ell_{1.5}$ & $\ell_2$ & $\ell_3$ \\\hline\cline{3-5}
\multirow{4}{*}{SUM} &  $\widehat{\bbeta}$  &  $(-0.2603, -0.9299, -1)$   &  $(-0.2603, -0.9299, -1)$   &  $(-0.2603, -0.9299, -1)$ \\
      %&  $\Phi^*$  &  $4854.747$  &  $5447.475$  &  $6110.581$\\
      &  $\GCoD$  &  $0.4133$  &  $0.3417$  &  $0.2615$\\
      &  $\%$  &  $62\%$  &  $62\%$  &  $62\%$\\
      &  $\epsilon_{90}$  &  $86.7791$  &  $86.7791$  &  $86.7791$\\ \hline
\multirow{4}{*}{MAX} &  $\widehat{\bbeta}$  &  $(-10.9038, -0.1571, -1)$   &  $(-10.9038, -0.1571, -1)$   &  $(10.9038, 0.1571, 1)$ \\
      %&  $\Phi^*$  &  $201.4506$  &  $207.1913$  &  $209.4644$\\
      &  $\GCoD$  &  $0.1821$  &  $0.1588$  &  $0.1495$\\
      &  $\%$  &  $10\%$  &  $10\%$  &  $10\%$\\
      &  $\epsilon_{90}$  &  $158.9295$  &  $158.9295$  &  $158.9295$\\ \hline
\multirow{4}{*}{SOS} &  $\widehat{\bbeta}$  &  $(2.4728, -0.2391, -1)$   &  $(-2.8551, 0.2102, 1)$   &  $(-3.1181, 0.1903, 1)$ \\
      %&  $\Phi^*$  &  $725637.8$  &  $790486.6$  &  $817725.3$\\
      &  $\GCoD$  &  $0.3163$  &  $0.2552$  &  $0.2295$\\
      &  $\%$  &  $8\%$  &  $8\%$  &  $8\%$\\
      &  $\epsilon_{90}$  &  $149.8204$  &  $151.9362$  &  $156.6873$\\ \hline
\multirow{4}{*}{1.5SUM} &  $\widehat{\bbeta}$  &  $(3.4138, -0.2225, -1)$   &  $(3.0670, -0.2704, -1)$   &  $(1.4864, -0.3260, -1)$ \\
      %&  $\Phi^*$  &  $73517.24$  &  $70879.39$  &  $64980.15$\\
      &  $\GCoD$  &  $0.1853$  &  $0.2145$  &  $0.2799$\\
      &  $\%$  &  $8\%$  &  $9\%$  &  $7\%$\\
      &  $\epsilon_{90}$  &  $149.6913$  &  $145.969$  &  $135.7776$\\ \hline
\multirow{4}{*}{kC} &  $\widehat{\bbeta}$  &  $(-2.6422, 0.2474, 1)$   &  $(-0.2632, -0.9011, -1)$   &  $(-0.3503, -0.9091, -1)$ \\
      %&  $\Phi^*$  &  $5671.552$  &  $5249.515$  &  $4679.65$\\
      &  $\GCoD$  &  $0.1263$  &  $0.1913$  &  $0.2791$\\
      &  $\%$  &  $9\%$  &  $57\%$  &  $60\%$\\
      &  $\epsilon_{90}$  &  $147.9623$  &  $84.4867$  &  $85.1974$\\ \hline
\multirow{4}{*}{AkC} &  $\widehat{\bbeta}$  &  $(-0.0741, 0.9357, 1)$   &  $(2.2028, 1.0126, 1)$   &  $(-0.9506, 0.9930, 1)$ \\
      %&  $\Phi^*$  &  $388.7408$  &  $309.9781$  &  $259.1419$\\
      &  $\GCoD$  &  $0.9468$  &  $0.9576$  &  $0.9645$\\
      &  $\%$  &  $64\%$  &  $70\%$  &  $65\%$\\
      &  $\epsilon_{90}$  &  $86.9840$  &  $94.2569$  &  $91.5147$\\ \hline
\multirow{4}{*}{MED} &  $\widehat{\bbeta}$  &  $(1.5779, -0.9545, -1)$   &  $(2.9207, 1.0139, 1)$   &  $(0.2899, 0.9792, 1)$ \\
      %&  $\Phi^*$  &  $6.865188$  &  $5.678452$  &  $5.046125$\\
      &  $\GCoD$  &  $0.9530$  &  $0.9611$  &  $0.9655$\\
      &  $\%$  &  $63\%$  &  $69\%$  &  $65\%$\\
      &  $\epsilon_{90}$  &  $88.5178$  &  $94.8548$  &  $90.5271$\\ \hline
\end{tabular}}
\end{table}

\begin{table}
\caption{Results for bidimensional experiments corrupting the $Y$ variables.\label{table:experiments6}}
\centering
{\scriptsize
\begin{tabular}{|c|c|c|c|c|}\cline{3-5}
       \multicolumn{2}{c|}{}            & V     & $\ell_1$ & $\ell_\infty$ \\\hline\cline{3-5}
\multirow{4}{*}{SUM} &  $\widehat{\bbeta}$  & $(-0.4324, -1.0070, -1)$   &  $(-2.7476, -1.1156, -1)$   &  $(-0.8817, -1.0333, -1)$   \\
      %&  $\Phi^*$  & $5764.311$  &  $5477.495$  &  $2853.437$  \\
      &  $\GCoD$  & $0.5226$  &  $0.5464$  &  $0.7637$  \\
      &  $\%$  & $72\%$  &  $57\%$  &  $73\%$  \\
      &  $\epsilon_{90}$  & $158.3495$  &  $144.4862$  &  $154.9621$  \\\hline
\multirow{4}{*}{MAX} &  $\widehat{\bbeta}$  & $(164.40, 1.95, -1)$   &  $(-131.52, -7.30, -1)$   &  $(-131.52, -7.30, -1)$   \\
      %&  $\Phi^*$  & $760.0988$  &  $186.3152$  &  $163.8737$  \\
      &  $\GCoD$  & $0.0109$  &  $0.7575$  &  $0.7867$  \\
      &  $\%$  & $5\%$  &  $6\%$  &  $6\%$  \\
      &  $\epsilon_{90}$  & $266.337$  &  $144.6019$  &  $144.6019$  \\\hline
\multirow{4}{*}{SOS} &  $\widehat{\bbeta}$  & $(-19.4780, 0.9765, 1)$   &  $(24.3778, -3.9704, -1)$   &  $(-21.8989, 2.4558, 1)$   \\
      %&  $\Phi^*$  & $2636508$  &  $679939$  &  $385975.8$  \\
      &  $\GCoD$  & $0.2459$  &  $0.8055$  &  $0.8896$  \\
      &  $\%$  & $24\%$  &  $12\%$  &  $14\%$  \\
      &  $\epsilon_{90}$  & $176.2108$  &  $119.0515$  &  $108.3728$  \\\hline
\multirow{4}{*}{1.5SUM} &  $\widehat{\bbeta}$  & $(2.2257, -0.9993, -1)$   &  $(8.1241, -2.8635, -1)$   &  $(4.2013, -1.5531, -1)$   \\
      %&  $\Phi^*$  & $114079$  &  $63837$  &  $35277.45$  \\
      &  $\GCoD$  & $0.3894$  &  $0.6583$  &  $0.8111$  \\
      &  $\%$  & $72\%$  &  $15\%$  &  $24\%$  \\
      &  $\epsilon_{90}$  & $161.1331$  &  $114.1084$  &  $107.9904$  \\\hline
\multirow{4}{*}{kC} &  $\widehat{\bbeta}$  & $(-0.6995, -0.9989, -1)$   &  $(4.8095, -1.6540, -1)$   &  $(-1.0107, -1.0744, -1)$   \\
      %&  $\Phi^*$  & $5599.338$  &  $5050.676$  &  $2745$  \\
      &  $\GCoD$  & $0.4422$  &  $0.4969$  &  $0.7265$  \\
      &  $\%$  & $71\%$  &  $23\%$  &  $67\%$  \\
      &  $\epsilon_{90}$  & $159.1129$  &  $100.6695$  &  $150.2014$  \\\hline
\multirow{4}{*}{AkC} &  $\widehat{\bbeta}$  & $(10.0084, -0.9838, -1)$   &  $(-1.3062, -1.0398, -1)$   &  $(-1.2815, -0.9942, -1)$   \\
      %&  $\Phi^*$  & $583.7418$  &  $286.736$  &  $129.5146$  \\
      &  $\GCoD$  & $0.7526$  &  $0.9914$  &  $0.9961$  \\
      &  $\%$  & $53\%$  &  $70\%$  &  $72\%$  \\
      &  $\epsilon_{90}$  & $168.5344$  &  $153.9189$  &  $159.2534$  \\\hline
\multirow{4}{*}{MED} &  $\widehat{\bbeta}$  & $(8.6545, -0.9641, -1)$   &  $(-0.8028, -1.0379, -1)$   &  $(-4.3252, -1.0113, -1)$   \\
      %&  $\Phi^*$  & $11.89961$  &  $7.084562$  &  $3.484066$  \\
      &  $\GCoD$  & $0.8478$  &  $0.9894$  &  $0.9947$  \\
      &  $\%$  & $57\%$  &  $73\%$  &  $69\%$  \\
      &  $\epsilon_{90}$  & $170.0131$  &  $154.4849$  &  $155.1026$  \\\hline
 \multicolumn{5}{c}{}\\\cline{3-5}
      \multicolumn{2}{c|}{}         & $\ell_{1.5}$ & $\ell_2$ & $\ell_3$ \\\hline\cline{3-5}
\multirow{4}{*}{SUM} &  $\widehat{\bbeta}$  &  $(-0.9890, -1.0403, -1)$   &  $(-0.9890, -1.0403, -1)$   &  $(-0.9890, -1.0403, -1)$ \\
      %&  $\Phi^*$  &  $4527.841$  &  $4034.632$  &  $3594.799$\\
      &  $\GCoD$  &  $0.6250$  &  $0.6658$  &  $0.7023$\\
      &  $\%$  &  $70\%$  &  $70\%$  &  $70\%$\\
      &  $\epsilon_{90}$  &  $154.0857$  &  $154.0857$  &  $154.0857$\\ \hline
\multirow{4}{*}{MAX} &  $\widehat{\bbeta}$  &  $(-131.52, -7.30, -1)$   &  $(-131.52, -7.30, -1)$   &  $(-131.52, -7.30, -1)$ \\
      %&  $\Phi^*$  &  $186.156$  &  $184.5921$  &  $180.275$\\
      &  $\GCoD$  &  $0.7577$  &  $0.7598$  &  $0.7654$\\
      &  $\%$  &  $6\%$  &  $6\%$  &  $6\%$\\
      &  $\epsilon_{90}$  &  $144.6019$  &  $144.6019$  &  $144.6019$\\ \hline
\multirow{4}{*}{SOS} &  $\widehat{\bbeta}$  &  $(24.0474, -3.7686, -1)$   &  $(23.2040, -3.2532, -1)$   &  $(22.5246, -2.8381, -1)$ \\
      %&  $\Phi^*$  &  $672228.4$  &  $631090.7$  &  $555242.8$\\
      &  $\GCoD$  &  $0.8077$  &  $0.8195$  &  $0.8412$\\
      &  $\%$  &  $13\%$  &  $13\%$  &  $13\%$\\
      &  $\epsilon_{90}$  &  $118.4519$  &  $119.827$  &  $115.0321$\\ \hline
\multirow{4}{*}{1.5SUM} &  $\widehat{\bbeta}$  &  $(8.2797, -2.4830, -1)$   &  $(5.8395, -1.9194, -1)$   &  $(4.7010, -1.6953, -1)$ \\
      %&  $\Phi^*$  &  $62263.44$  &  $56482.47$  &  $48865.13$\\
      &  $\GCoD$  &  $0.6667$  &  $0.6976$  &  $0.7384$\\
      &  $\%$  &  $14\%$  &  $19\%$  &  $23\%$\\
      &  $\epsilon_{90}$  &  $114.0191$  &  $102.4955$  &  $97.65193$\\ \hline
\multirow{4}{*}{kC} &  $\widehat{\bbeta}$  &  $(-1.0107, -1.0744, -1)$   &  $(-1.0107, -1.0744, -1)$   &  $(-0.8903, -1.0744, -1)$ \\
      %&  $\Phi^*$  &  $4351.828$  &  $3879.523$  &  $3457.372$\\
      &  $\GCoD$  &  $0.5665$  &  $0.6135$  &  $0.6556$\\
      &  $\%$  &  $67\%$  &  $67\%$  &  $66\%$\\
      &  $\epsilon_{90}$  &  $150.2014$  &  $150.2014$  &  $150.2834$\\ \hline
\multirow{4}{*}{AkC} &  $\widehat{\bbeta}$  &  $(-2.6754, -1.0658, -1)$   &  $(-2.7011, -0.9640, -1)$   &  $(-3.9149, -1.0070, -1)$ \\
      %&  $\Phi^*$  &  $329.7549$  &  $300.3972$  &  $282.1704$\\
      &  $\GCoD$  &  $0.9901$  &  $0.9910$  &  $0.9915$\\
      &  $\%$  &  $69\%$  &  $68\%$  &  $69\%$\\
      &  $\epsilon_{90}$  &  $150.0206$  &  $161.8515$  &  $155.8964$\\ \hline
\multirow{4}{*}{MED} &  $\widehat{\bbeta}$  &  $(-0.8019, -1.0319, -1)$   &  $(-2.6799, -1.0009, -1)$   &  $(-1.5141, -1.0345, -1)$ \\
      %&  $\Phi^*$  &  $5.904683$  &  $5.028813$  &  $4.785246$\\
      &  $\GCoD$  &  $0.9911$  &  $0.9924$  &  $0.9928$\\
      &  $\%$  &  $74\%$  &  $70\%$  &  $70\%$\\
      &  $\epsilon_{90}$  &  $155.184$  &  $157.4707$  &  $154.3846$\\ \hline
\end{tabular}}
\end{table}

\begin{table}
\caption{Results for Experiments for $d=4$ and corrupting the $X$ variables.\label{table:experiments4}}
\centering
{\tiny
%\hspace*{-1.75cm}
\begin{tabular}{|c|c|c|c|c|}\cline{3-5}
           \multicolumn{2}{c|}{}        & V     & $\ell_1$ & $\ell_\infty$ \\\hline\cline{3-5}
\multirow{4}{*}{SUM} &  $\widehat{\bbeta}$  & $(8.7754, 0.2361, 0.1242, -0.0645, 1)$   &  $(-167.9861, 32.8678, -11.1472, -15.3593, 1)$   &  $(19.6624, 1.9411, 1.4336, -2.6949, 1)$   \\
      %&  $\Phi^*$  & $13499.58$  &  $9073.616$  &  $4162.037$  \\
      &  $\GCoD$  & $0.0369$  &  $0.3527$  &  $0.7030$  \\
      &  $\%$  & $8\%$  &  $9\%$  &  $15\%$  \\
      &  $\epsilon_{90}$  & $285.1339$  &  $172.616$  &  $166.2396$  \\\hline
\multirow{4}{*}{MAX} &  $\widehat{\bbeta}$  & $(11.2676, -0.8055, 0.4093, 0.3802, 1)$   &  $(95.4943, -2.3074, -2.7088, 4.5984, 1)$   &  $(76.9688, -2.1455, -2.9597, 4.6480, 1)$   \\
      %&  $\Phi^*$  & $401.0232$  &  $226.1476$  &  $97.87958$  \\
      &  $\GCoD$  & $0.1200$  &  $0.5037$  &  $0.7852$  \\
      &  $\%$  & $2\%$  &  $9\%$  &  $6\%$  \\
      &  $\epsilon_{90}$  & $243.9038$  &  $160.86$  &  $164.3572$  \\\hline
\multirow{4}{*}{SOS} &  $\widehat{\bbeta}$  & $(2.7637, 0.1306, 0.06391, -0.0111, 1)$   &  $(-35.0079, -17.4180, 5.1138, 8.8243, -1)$   &  $(14.4492, 2.3985, 1.8254, -3.4712, 1)$   \\
      %&  $\Phi^*$  & $3017114$  &  $1325324$  &  $287765.5$  \\
      &  $\GCoD$  & $0.0409$  &  $0.5787$  &  $0.9085$  \\
      &  $\%$  & $6\%$  &  $9\%$  &  $8\%$  \\
      &  $\epsilon_{90}$  & $285.0815$  &  $170.37$  &  $165.6255$  \\\hline
\multirow{4}{*}{1.5SUM} &  $\widehat{\bbeta}$  & $(3.1382, 0.1714, 0.0663, -0.0352 1)$   &  $(21.9152, -18.9245, 5.5144, 9.6284, -1)$   &  $(-20.1562, -2.0728, -1.5407, 2.9444, -1)$   \\
      %&  $\Phi^*$  & $193942.4$  &  $105739.7$  &  $33409.83$  \\
      &  $\GCoD$  & $0.0418$  &  $0.4776$  &  $0.8349$  \\
      &  $\%$  & $7\%$  &  $8\%$  &  $14\%$  \\
      &  $\epsilon_{90}$  & $282.7383$  &  $167.7096$  &  $165.9725$  \\\hline
\multirow{4}{*}{kC} &  $\widehat{\bbeta}$  & $(-6.8937, 0.1108, 0.0744, -0.0183, 1)$   &  $(-34.1432, -15.4977, 4.3066, 7.9523, -1)$   &  $(5.0421, 2.0898, 1.4381, -2.8638, 1)$   \\
      %&  $\Phi^*$  & $10943.27$  &  $7315.619$  &  $3387.323$  \\
      &  $\GCoD$  & $0.0258$  &  $0.3487$  &  $0.6984$  \\
      &  $\%$  & $8\%$  &  $8\%$  &  $15\%$  \\
      &  $\epsilon_{90}$  & $276.4327$  &  $168.3023$  &  $169.65$  \\\hline
\multirow{4}{*}{AkC} &  $\widehat{\bbeta}$  & $(-29.5486, 0.5489, 0.2119, 0.2342, 1)$   &  $(11.5813, 2.8055, -0.1579, 0.1805, 1)$   &  $(2.7269, 1.0225, 0.9985, 1.0072, 1)$   \\
      %&  $\Phi^*$  & $2665.219$  &  $2283.138$  &  $88.06644$  \\
      &  $\GCoD$  & $0.1544$  &  $0.8716$  &  $0.9950$  \\
      &  $\%$  & $12\%$  &  $5\%$  &  $82\%$  \\
      &  $\epsilon_{90}$  & $304.1316$  &  $306.9669$  &  $496.6216$  \\\hline
\multirow{4}{*}{MED} &  $\widehat{\bbeta}$  & $(11.3163, 0.5095, 0.5018, 0.0667, 1)$   &  $(15.2913, -1.38181, -0.1062, 9.6624, 1)$   &  $(2.3001, 1.0447, 1.0149, 1.0033, 1)$   \\
      %&  $\Phi^*$  & $68.29634$  &  $60.18373$  &  $2.075812$  \\
      &  $\GCoD$  & $0.3706$  &  $0.8308$  &  $0.9941$  \\
      &  $\%$  & $9\%$  &  $11\%$  &  $80\%$  \\
      &  $\epsilon_{90}$  & $283.331$  &  $251.5948$  &  $497.3323$  \\\hline
       \multicolumn{5}{c}{}\\ \cline{3-5}
    \multicolumn{2}{c|}{}      & $\ell_{1.5}$ & $\ell_2$ & $\ell_3$ \\\hline\cline{3-5}
\multirow{4}{*}{SUM} &  $\widehat{\bbeta}$  &  $(-25.3339, 7.2803, 0.3850, -6.5208, 1)$   &  $(-25.3339, 7.2803, 0.3850, -6.5208, 1)$   &  $(-48.9741, -2.5251, -1.5173, 3.4889, -1)$ \\
      %&  $\Phi^*$  &  $8448.604$  &  $7526.85$  &  $6383.166$\\
      &  $\GCoD$  &  $0.3973$  &  $0.4630$  &  $0.5446$\\
      &  $\%$  &  $12\%$  &  $12\%$  &  $11\%$\\
      &  $\epsilon_{90}$  &  $167.1534$  &  $167.1534$  &  $163.8287$\\ \hline
\multirow{4}{*}{MAX} &  $\widehat{\bbeta}$  &  $(-76.9688, 2.1455, 2.9597, -4.6480, -1)$   &  $(-76.9688, 2.1455, 2.9597, -4.6480, -1)$   &  $(-76.9688, 2.1455, 2.9597, -4.6480, -1)$ \\
      %&  $\Phi^*$  &  $204.599$  &  $177.8851$  &  $147.1322$\\
      &  $\GCoD$  &  $0.5510345$  &  $0.6096547$  &  $0.677138$\\
      &  $\%$  &  $6\%$  &  $6\%$  &  $6\%$\\
      &  $\epsilon_{90}$  &  $164.3572$  &  $164.3572$  &  $164.3572$\\ \hline
\multirow{4}{*}{SOS} &  $\widehat{\bbeta}$  &  $(-19.8365, -24.1780, -1.6843, 23.0309, -1)$   &  $(-37.1798, -20.6518, -4.8914, 22.4924, -1)$   &  $(16.2930, 4.1351, 2.2042, -5.3890, 1)$ \\
      %&  $\Phi^*$  &  $1135243$  &  $896820.1$  &  $653754.7$\\
      &  $\GCoD$  &  $0.6391$  &  $0.7149$  &  $0.7921$\\
      &  $\%$  &  $9\%$  &  $9\%$  &  $4\%$\\
      &  $\epsilon_{90}$  &  $159.013$  &  $160.1321$  &  $165.3201$\\ \hline
\multirow{4}{*}{1.5SUM} &  $\widehat{\bbeta}$  &  $(27.4692, 14.0582, 1.0081, -12.9659, 1)$   &  $(27.4555, 14.0608, 1.0082, -12.9683, 1)$   &  $(-20.4048, -3.2308, -1.6763, 4.1796, -1)$ \\
      %&  $\Phi^*$  &  $94853.7$  &  $79761.93$  &  $62563.54$\\
      &  $\GCoD$  &  $0.5314$  &  $0.6059$  &  $0.6909$\\
      &  $\%$  &  $10\%$  &  $10\%$  &  $5\%$\\
      &  $\epsilon_{90}$  &  $162.8882$  &  $162.8875$  &  $164.1443$\\ \hline
\multirow{4}{*}{kC} &  $\widehat{\bbeta}$  &  $(31.8219, 41.5015, -5.2288, -30.4070, 1)$   &  $(2.4227, 14.3655, 4.4768, -15.4827, 1)$   &  $(6.6713, -3.7849, -1.5627, 4.3751, -1)$ \\
      %&  $\Phi^*$  &  $6834.039$  &  $6032.707$  &  $5121.862$\\
      &  $\GCoD$  &  $0.3916$  &  $0.4629$  &  $0.5440$\\
      &  $\%$  &  $5\%$  &  $7\%$  &  $4\%$\\
      &  $\epsilon_{90}$  &  $165.793$  &  $168.1855$  &  $165.9668$\\ \hline
\multirow{4}{*}{AkC} &  $\widehat{\bbeta}$  &  $(7.9530, -1.6065, 0.3482, 0.8960, -1)$   &  $(-25.2618, -1.0371, -1.4553, 0.7368, -1)$   &  $(40.7617, -1.6662, -0.5106, 0.5624, -1)$ \\
      %&  $\Phi^*$  &  $4618.866$  &  $3292.832$  &  $2103.07$\\
      &  $\GCoD$  &  $0.7403$  &  $0.8148$  &  $0.8817$\\
      &  $\%$  &  $7\%$  &  $11\%$  &  $9\%$\\
      &  $\epsilon_{90}$  &  $180.9401$  &  $244.0442$  &  $231.9954$\\ \hline
\multirow{4}{*}{MED} &  $\widehat{\bbeta}$  &  $(-28.1536, -1.9062, -0.5785, 0.5246, -1)$   &  $(-51.5261, 1.9897, 1.0285, -0.5282, 1)$   &  $(6.9522, 1.2873, 1.0511, -0.1044, 1)$ \\
      %&  $\Phi^*$  &  $61.24124$  &  $50.68379$  &  $37.6414$\\
      &  $\GCoD$  &  $0.8278$  &  $0.8575$  &  $0.8941$\\
      &  $\%$  &  $9\%$  &  $8\%$  &  $14\%$\\
      &  $\epsilon_{90}$  &  $237.8898$  &  $305.539$  &  $350.0691$\\ \hline
\end{tabular}}
\end{table}

\begin{table}
\caption{Results for Experiments for $d=4$ and corrupting the $Y$ variables.\label{table:experiments8}}
\centering
{\tiny
\hspace*{-1.75cm}
\begin{tabular}{|c|c|c|c|c|}\cline{3-5}
       \multicolumn{2}{c|}{}            & V     & $\ell_1$ & $\ell_\infty$ \\\hline\cline{3-5}
\multirow{4}{*}{SUM} &  $\widehat{\bbeta}$  & $(1.9468, 0.9648, 0.9899, 1.0058, 1)$   &  $(-1.9158, -1.1083, -0.8751, -3.3186, -1)$   &  $(1.6655, -1.0083, -1.0530, -1.0446, -1)$   \\
      %&  $\Phi^*$  & $7653.35$  &  $6622.243$  &  $1900.56$  \\
      &  $\GCoD$  & $0.5999$  &  $0.6538$  &  $0.9006$  \\
      &  $\%$  & $78\%$  &  $14\%$  &  $76\%$  \\
      &  $\epsilon_{90}$  & $123.5456$  &  $149.6274$  &  $121.8106$  \\\hline
\multirow{4}{*}{MAX} &  $\widehat{\bbeta}$  & $(1-04.7766, -1.0780, -2.8506, -0.8355, -1)$   &  $(120.6153, -1.4207, -5.5268, -0.7782, -1)$   &  $(54.3395, 2.3207, 6.0411, 3.4977, 1)$   \\
      %&  $\Phi^*$  & $703.3325$  &  $183.4559$  &  $97.5249$  \\
      &  $\GCoD$  & $0.3357$  &  $0.8267$  &  $0.9078$  \\
      &  $\%$  & $12\%$  &  $7\%$  &  $12\%$  \\
      &  $\epsilon_{90}$  & $151.6067$  &  $147.4952$  &  $138.4277$  \\\hline
\multirow{4}{*}{SOS} &  $\widehat{\bbeta}$  & $(-12.1432, -0.8507, -1.0758, -1.1049, -1)$   &  $(25.1165, -1.2149, -5.4326, -1.1199, -1)$   &  $(-5.4787, -1.8048, -2.3397, -2.0389, -1)$   \\
      %&  $\Phi^*$  & $4347478$  &  $743813.1$  &  $150116.9$  \\
      &  $\GCoD$  & $0.4247$  &  $0.9015$  &  $0.9801$  \\
      &  $\%$  & $45\%$  &  $13\%$  &  $15\%$  \\
      &  $\epsilon_{90}$  & $124.0456$  &  $135.9287$  &  $102.1587$  \\\hline
\multirow{4}{*}{1.5SUM} &  $\widehat{\bbeta}$  & $(-2.1265, -0.9557, -0.9984, -1.0235, -1)$   &  $(34.3751, -1.0783, -5.2458, -1.0619, -1)$   &  $(-0.6651, -1.3869, -1.5549, -1.5790, -1)$   \\
      %&  $\Phi^*$  & $172667.9$  &  $68991.48$  &  $18141.06$  \\
      &  $\GCoD$  & $0.5106$  &  $0.8044$  &  $0.9485$  \\
      &  $\%$  & $77\%$  &  $11\%$  &  $22\%$  \\
      &  $\epsilon_{90}$  & $124.3694$  &  $139.4734$  &  $95.54551$  \\\hline
\multirow{4}{*}{kC} &  $\widehat{\bbeta}$  & $(-0.3095, -0.9816, -1.0017, -1.009643, -1)$   &  $(2.1980, -0.8680, -0.9950, -3.4086, -1)$   &  $(-0.6929, -1.0211, -1.0606, -1.0666, -1)$   \\
      %&  $\Phi^*$  & $7469.548$  &  $5493.28$  &  $1841.433$  \\
      &  $\GCoD$  & $0.5275$  &  $0.6525$  &  $0.8835$  \\
      &  $\%$  & $80\%$  &  $10\%$  &  $74\%$  \\
      &  $\epsilon_{90}$  & $123.0891$  &  $145.6142$  &  $120.8033$  \\\hline
\multirow{4}{*}{AkC} &  $\widehat{\bbeta}$  & $(-7.2126, -0.9981, -1.2345, -0.9988, -1)$   &  $(-1.7307, -0.9801, -1.0396, -1.0121, -1)$   &  $(0.1128, -0.9847, -1.0149, -1.0013, -1)$   \\
      %&  $\Phi^*$  & $700.2851$  &  $321.0255$  &  $90.39922$  \\
      &  $\GCoD$  & $0.8785$  &  $0.9933$  &  $0.9981$  \\
      &  $\%$  & $57\%$  &  $77\%$  &  $80\%$  \\
      &  $\epsilon_{90}$  & $105.7586$  &  $120.4785$  &  $121.9634$  \\\hline
\multirow{4}{*}{MED} &  $\widehat{\bbeta}$  & $(-8.4437, -1.0328, -1.1891, -0.9958, -1)$   &  $(-3.0605, -0.9660 -1.0175, -1.0366, -1)$   &  $(-1.7471, -0.9713, -0.9881, -1.0144, -1)$   \\
      %&  $\Phi^*$  & $15.18326$  &  $7.56219$  &  $1.897255$  \\
      &  $\GCoD$  & $0.9011$  &  $0.9921$  &  $0.9980$  \\
      &  $\%$  & $58\%$  &  $76\%$  &  $79\%$  \\
      &  $\epsilon_{90}$  & $105.9371$  &  $123.0289$  &  $123.8959$  \\\hline
 \multicolumn{5}{c}{}\\\cline{3-5}
      \multicolumn{2}{c|}{}         & $\ell_{1.5}$ & $\ell_2$ & $\ell_3$ \\\hline\cline{3-5}
\multirow{4}{*}{SUM} &  $\widehat{\bbeta}$  &  $(0.5934, -1.0202, -1.0588, -1.0264, -1)$   &  $(0.6616, -1.0203, -1.0584, -1.0270, -1)$   &  $(0.9775, -1.0098, -1.0563, -1.0343, -1)$ \\
      %&  $\Phi^*$  &  $4802.058$  &  $3813.333$  &  $3025.717$\\
      &  $\GCoD$  &  $0.7489$  &  $0.8006$  &  $0.8418$\\
      &  $\%$  &  $80\%$  &  $80\%$  &  $78\%$\\
      &  $\epsilon_{90}$  &  $119.4431$  &  $119.5293$  &  $120.6788$\\ \hline
\multirow{4}{*}{MAX} &  $\widehat{\bbeta}$  &  $(120.6153, -1.4207, -5.5268, -0.7782, -1)$   &  $(-54.3395, -2.3207, -6.0411, -3.4977, -1)$   &  $(-54.3395, -2.3207, -6.0411, -3.4977, -1)$ \\
      %&  $\Phi^*$  &  $183.4559$  &  $171.0524$  &  $143.6565$\\
      &  $\GCoD$  &  $0.8267$  &  $0.8384$  &  $0.8643$\\
      &  $\%$  &  $7\%$  &  $12\%$  &  $12\%$\\
      &  $\epsilon_{90}$  &  $147.4952$  &  $138.4277$  &  $138.4277$\\ \hline
\multirow{4}{*}{SOS} &  $\widehat{\bbeta}$  &  $(-14.4853, 1.5436, 4.4201, 1.5950, 1)$   &  $(-0.3904, 1.7361, 2.9264, 2.0617, 1)$   &  $(4.7620, 1.9721, 2.5444, 2.0415, 1)$ \\
      %&  $\Phi^*$  &  $738909.1$  &  $549957.6$  &  $367130$\\
      &  $\GCoD$  &  $0.9022$  &  $0.9272$  &  $0.9514$\\
      &  $\%$  &  $13\%$  &  $10\%$  &  $12\%$\\
      &  $\epsilon_{90}$  &  $131.3351$  &  $114.7621$  &  $106.4697$\\ \hline
\multirow{4}{*}{1.5SUM} &  $\widehat{\bbeta}$  &  $(15.7120, -1.1641, -2.6186, -1.8366, -1)$   &  $(-0.8627, -1.4497, -1.6239, -1.9098, -1)$   &  $(-0.6434, -1.4056, -1.5798, -1.5348, -1)$ \\
      %&  $\Phi^*$  &  $67759.55$  &  $50632.87$  &  $36514.17$\\
      &  $\GCoD$  &  $0.8079$  &  $0.8565$  &  $0.8965$\\
      &  $\%$  &  $21\%$  &  $22\%$  &  $20\%$\\
      &  $\epsilon_{90}$  &  $114.939$  &  $97.67539$  &  $97.29497$\\ \hline
\multirow{4}{*}{kC} &  $\widehat{\bbeta}$  &  $(-1.0976, -1.0234, -1.0643, -1.0656, -1)$   &  $(-1.0942, -1.0234, -1.0641, -1.0656, -1)$   &  $(-0.7613, -1.0216, -1.0617, -1.0665, -1)$ \\
      %&  $\Phi^*$  &  $4657.853$  &  $3697.56$  &  $2933.1$\\
      &  $\GCoD$  &  $0.7053$  &  $0.7661$  &  $0.8144$\\
      &  $\%$  &  $74\%$  &  $74\%$  &  $74\%$\\
      &  $\epsilon_{90}$  &  $120.25$  &  $120.262$  &  $120.6901$\\ \hline
\multirow{4}{*}{AkC} &  $\widehat{\bbeta}$  &  $(0.8072, -0.9319, -1.1111, -1.0901, -1)$   &  $(-1.5573, -0.9672, -0.9991, -1.0184, -1)$   &  $(2.4443, -1.0165, -0.9923, -1.0147, -1)$ \\
      %&  $\Phi^*$  &  $335.7945$  &  $219.6891$  &  $330.8805$\\
      &  $\GCoD$  &  $0.9929$  &  $0.9954$  &  $0.9930$\\
      &  $\%$  &  $64\%$  &  $77\%$  &  $82\%$\\
      &  $\epsilon_{90}$  &  $124.0139$  &  $123.7847$  &  $123.5452$\\ \hline
\multirow{4}{*}{MED} &  $\widehat{\bbeta}$  &  $(-0.6735, -0.9887, -1.0180, -0.9497, -1)$   &  $(0.4156, -0.9995, -1.0147, -1.0116, -1)$   &  $(-1.1572, -0.9753, -1.0309, -0.9853, -1)$ \\
      %&  $\Phi^*$  &  $5.2003$  &  $4.875409$  &  $3.435911$\\
      &  $\GCoD$  &  $0.9945$  &  $0.9949$  &  $0.9964$\\
      &  $\%$  &  $75\%$  &  $81\%$  &  $78\%$\\
      &  $\epsilon_{90}$  &  $118.3319$  &  $121.9701$  &  $120.0091$\\ \hline
\end{tabular}}
\end{table}
%\end{landscape} 

Tables \ref{table:experiments2}-\ref{table:experiments8} report, for each  battery of generated data, the following information: i) the coefficients of the optimal hyperplane ($\widehat{\bbeta}$), ii) the goodness of fitting index $\GCoD$, iii) the percentage of the sample data  which are contained in a strip delimited by two parallel hyperplanes to $y =\widehat{\bbeta}x$  with  (orthogonal) distance $\varepsilon=10$  ($\%$), and iv) the width of the strip that is necessary to include $90\%$ of the data ($\epsilon_{90}$).

We conclude, from the experiments for the bivariate case, that in general a better performance is observed in all the methods when the corrupted  coordinate is the dependent one ($Y$), as compared with introducing the corruption on the independent  coordinate ($X$). In particular, the SUM, the 1.5SUM and the kC criteria (for vertical distance residuals) get better fitting models in the $Y$-corrupted case. Although slightly better, almost similar results were obtained for the AkC, MEDIAN and kC (for $\ell_\tau$ residuals) due to the robustness of those criteria. Also, we observe that for the $X$-corrupted  case, the linear residuals (V, $\ell_1$ and $\ell_\infty$) models coincide for all the criteria except the AkC. This is not the case in the $Y$-corrupted  experiments, where equal or similar models were obtained for all the $\ell_\tau$-residuals. Observe that although in the $X$-corrupted  case the larger $\%$ seems to imply a greater $\GCoD$, that is not the case in the $Y$-corrupted experiments where one can find many combinations of criteria-residuals where that behavior does not happen.

Similar conclusions can be derived from the multivariate case ($d=4$), except that in this case there are no coincidences between the models obtained with different combinations of criteria and residuals. Furthermore, the convenience of using measures for the goodness of fitting which are not criterion/residual dependent is confirmed.

\subsection{Data: Durbin-Watson}

We also performed some experiments over the classical real data sample used in \cite{dw}. The data aims to analyze the annual consumption of spirits from 1870 to 1938 ($n=69$) from the incomes and the relative price of spirits (deflated by a cost-of-living index). Hence, the variables observed in this data sets are the logarithms (the coefficients are then interpreted in terms of percent change) of the following measures: $X_1$ (Real income per head),  $X_2$ (Relative price of spirits) and $X_3$ (Consumption of spirits per head).

For illustrative purposes, we analyze both the global model with the three variables ($d=3$) and the bivariate model considering $X_1$ and $X_3$ and obviating $X_2$ ($d=2$).

\subsubsection{Bivariate model}
First, for the case $d=2$, we run the 42  models (Table \ref{tablef}) over the data set where $X_1$ (income) and $X_3$ (consumption) are measured. The obtained hyperplanes are detailed in Table \ref{table:dw1} and the fitted lines drawn in Figure \ref{fig:dw1}. Note that the methods that use vertical distance residuals were not able to capture the actual behavior of the consumption with respect to the incomes. Furthermore, the MAX criterion seems to fail for any choice of  residuals, since it tries to explain the unique outlier point that exists in the data set. The rest of the hyperplanes, with minimal deviations, have a similar behavior.  In order to analyze the differences between these models we also report in Table \ref{table:mv} the marginal variations of each one of the models (according to Lemma \ref{le:response}).

\input{figura-dw1.tex}

Observe that, when the $\ell_1$ residuals are considered, all except the MAX criterion provide a $0$ marginal variation. This pattern can be explained as a result of Lemma \ref{le:responsesnorms} and the fact that the $\ell_1$-norm unit ball in $\R^2$ has extreme points $\{\pm(0,1), \pm(1,0)\}$. Hence $\k(\beta) = \left\{\begin{array}{cl}
1 & \mbox{if $\beta_2=\max\{|\beta_1|, |\beta_2|\}$},\\
-1 & \mbox{if $\beta_2=-\max\{|\beta_1|, |\beta_2|\}$},\\
0 & \mbox{otherwise.}
\end{array}\right.$. Thus, the marginal variation of $X_1$ with respect to $X_3$ is zero iff $|\beta_1| = \max\{|\beta_1|, |\beta_2|\}$, being then $|\beta_2| < |\beta_1|$. It means that the absolute value of the slope of the line is greater than $1$, being the decreasing (or increasing) of the response consumption in terms of the incomes more than a $100\%$.

In order to validate and analyze the stability of the computed hyperplanes we perform a $k$-fold cross validation scheme \cite{kfold} to the data set. Such a method consists of randomly partitioning the sample into $k$ folds of similar size, $S_1, \ldots, S_k$. For each $j \in \{1, \ldots, k\}$, each optimal hyperplane is computed using the points in $\bigcup_{i\neq j} S_i$ and $S_j$ is used to validate the results. In our case, we partitioned the data into $k=7$ folds, each of them with $10$ data, except one with $9$ points. In Table \ref{table:dw2} we summarize the results obtained with this experiment. We report: the maximum, minimum,  median and mean width of the strips that are necessary to cover the  $90\%$ of the (validation) data for the seven runs.

\input{figura-dw2.tex}

From the above results, we note that the models that use vertical distance residuals need, in general larger strips to cover the $90\%$ of the points. The strips are particulary large for the MEDIAN criterion, where the widest strips were obtained. This conclusion is justified since the quantile criteria accommodate a single point, but do not take into account the deviations to the remainder elements in the data (apart from the ordering in the residuals). Also, for the same reason, the conservative MAX criterion makes the models to require wider strips. The main observed difference between the MEDIAN and the MAX criteria is that whereas the behavior (in term of the fitting strips) of the MAX criterion is similar for the six choices of residuals, the MEDIAN gets very different results depending of the chosen residual. The most robust residuals, based on the smallest range between the maximum and minimum length of the strips, are the $\ell_{1}$,  $\ell_{1.5}$, and $\ell_3$; while with the same measure of robustness, the $k$-centrum criterion gets the best results.

To illustrate the quality of the optimal hyperplanes, in Figure \ref{fig:dw2} we show the values of the consumptions with respect to the actual consumptions for the first random fold in the experiments (in the validation sample that was not used to compute the hyperplanes).

The conclusions are that the vertical distance residuals do not fit well to the actual the trend of the validation data. The same conclusion also applies to the models that use the MAX criterion or $\ell_\infty$ residuals. On the other hand, the $\ell_1$-residual models seem to fit quite well to the data, whereas the $\ell_\tau$-residual models have similar (good) behavior. As expected the $kC$ and $AkC$ criteria, which are known to be very robust, actually capture the main information about the trend of the data.

\subsubsection{Complete models}

We also performed the same experiments for the whole data set. The three variables $X_1$ (incomes), $X_2$ (prices) and $X_3$ (consumptions) are now considered. The optimal hyperplanes are shown in Table \ref{table:dw3} (since the coefficients are non zero they were divided by $-\beta_3$ to make easier the interpretation and representations of the models as $X_3=\beta_0 + \beta_1X_1+\beta_2X_2$).

The summary of the results of the $k$-fold cross validation scheme (where the data set was partitioned exactly as in the bivariate case) is shown in Table \ref{table:dw4}.  Finally, Figure \ref{fig:dw3} shows the values of the consumptions with respect to the actual consumptions for the first random fold in the experiments. From the results, one can observe that including all the variables in the model reduces the differences among the models obtained with  the different methods. In this case, the consumption seems to be well linearly described by the incomes and prices. This conclusion is supported both by the projection and by the summary of k-cross validation experiments. The exceptionally bad performance of the MAX criterion in the former case (the model that only included $X_1$ and $X_3$), is now as good as the rest of the criteria. In addition, the inclusion of prices in the model fixes the, in most cases, senseless  signs of the coefficients in the simple models in Table \ref{table:mv}. One can observe that in those cases an increase of the incomes would predict a decrease of the consumptions. This unusual trend is fixed by introducing the prices in the complete model.

\begin{table}
\caption{Estimations for the Durbin-Watson's dataset.\label{table:dw3}}
\centering
{\tiny
\hspace*{-0.2cm}\begin{tabular}{|c|c|c|c|}\cline{2-4}
   \multicolumn{1}{c|}{}   &  V    &  $\ell_1$  &  $\ell_\infty$   \\\hline
SUM   & $(4.4817, 0.0696, -1.3374, -1)$ & $(4.555, 0.0587, -1.3623, -1)$ & $(4.1367, 0.3502, -1.4305, -1)$ \\
MAX   & $(4.5227, 0.0646, -1.3519, -1)$ & $(4.6159, -0.013, -1.3273, -1)$ & $(4.1355, 0.5086, -1.5758, -1)$ \\
SOS   & $(3.9725, 0.0331, -1.0692, -1)$ & $(4.404, 0.1369, -1.3881, -1)$ & $(4.404, 0.1369, -1.3881, -1)$ \\
1.5SUM & $(4.404, 0.1369, -1.3881, -1)$ & $(4.404, 0.1369, -1.3881, -1)$ & $(4.404, 0.1369, -1.3881, -1)$ \\
kC    & $(4.4159, 0.0288, -1.2753, -1)$ & $(4.4905, 0.0635, -1.3425, -1)$ & $(4.3334, 0.1325, -1.3317, -1)$ \\
AkC   & $(4.4355, 0.0655, -1.3183, -1)$ & $(4.4521, 0.0585, -1.3197, -1)$ & $(4.4688, 0.0535, -1.323, -1)$ \\
MED & $(4.4288, 0.0488, -1.2979, -1)$ & $(4.5075, 0.0634, -1.3476, -1)$ & $(4.3559, 0.1431, -1.3489, -1)$ \\\hline
 \multicolumn{4}{c}{}\\\cline{2-4}
   \multicolumn{1}{c|}{} &  $\ell_{1.5}$   &  $\ell_{2}$   &   $\ell_{3}$  \\\hline
SUM   & $(4.4445, 0.0698, -1.3242, -1)$ & $(4.472, 0.0633, -1.331, -1)$ & $(4.4922, 0.0619, -1.3386, -1)$ \\
MAX   & $(4.4155, 0.0352, -1.2797, -1)$ & $(4.3938, 0.1107, -1.3377, -1)$ & $(4.2655, 0.1691, -1.3326, -1)$ \\
SOS   & $(4.3498, 0.1131, -1.3201, -1)$ & $(4.3498, 0.1131, -1.3201, -1)$ & $(4.3498, 0.1131, -1.3201, -1)$ \\
1.5SUM & $(4.2123, 0.4308, -1.5386, -1)$ & $(4.0853, 0.4429, -1.4891, -1)$ & $(3.6048, 0.7761, -1.5744, -1)$ \\
kC    & $(5.2647, -0.6758, -1.0312, -1)$ & $(3.5719, 1.1094, -1.8642, -1)$ & $(3.4912, 1.0623, -1.7796, -1)$ \\
AkC   & $(4.1061, 0.5015, -1.551, -1)$ & $(4.1579, 0.467, -1.5434, -1)$ & $(4.2963, 0.3239, -1.4761, -1)$ \\
MED & $(4.3576, 0.2689, -1.4559, -1)$ & $(4.0772, 0.4066, -1.4415, -1)$ & $(76.3635, 25.0913, -61.4268, -1)$ \\\hline
\end{tabular}}
\end{table}

\setlength\tabcolsep{4pt}

\begin{table}
\caption{Summary of k-fold cross validations experiments for the Durbin-Watson's dataset.\label{table:dw4}}
\centering
{\tiny\begin{tabular}{|c|c|c|c|c|c|c|c|}\cline{3-8}
\multicolumn{2}{c|}{}  & V & $\ell_1$ & $\ell_\infty$  & $\ell_{1.5}$  & $\ell_{2}$  &  $\ell_{3}$ \\\hline
 \multirow{4}{*}{SUM} & $\min \varepsilon_{90}$ & 0.0369 & 0.0388 & 0.0315 & 0.0380 & 0.0346 & 0.0347 \\
      & $\max \varepsilon_{90}$ & 0.0735 & 0.0741 & 0.0832 & 0.0743 & 0.0743 & 0.0732 \\
      & ${\rm median} \varepsilon_{90}$ & 0.0629 & 0.0627 & 0.0647 & 0.0625 & 0.0625 & 0.0626 \\
      & $\varepsilon_{90}$ & 0.0573 & 0.0598 & 0.0616 & 0.0580 & 0.0567 & 0.0593 \\
\hline
\multirow{4}{*}{MAX} & $\min \varepsilon_{90}$ & 0.0562 & 0.0515 & 0.0515 & 0.0515 & 0.0515 & 0.0515 \\
      & $\max \varepsilon_{90}$ & 0.0807 & 0.0762 & 0.0760 & 0.0760 & 0.0760 & 0.0762 \\
      & ${\rm median} \varepsilon_{90}$ & 0.0701 & 0.0607 & 0.0644 & 0.0644 & 0.0607 & 0.0607 \\
      & $\varepsilon_{90}$ & 0.0678 & 0.0624 & 0.0641 & 0.0641 & 0.0624 & 0.0624 \\
\hline
\multirow{4}{*}{SOS} & $\min \varepsilon_{90}$ & 0.0255 & 0.0362 & 0.0310 & 0.0321 & 0.0327 &0.0327 \\
      & $\max \varepsilon_{90}$ & 0.0656 & 0.0683 & 0.0691 & 0.0678 & 0.0675 & 0.0675 \\
      & ${\rm median} \varepsilon_{90}$ & 0.0586 & 0.0583 & 0.0568 & 0.0586 & 0.0581 & 0.0582 \\
      & $\varepsilon_{90}$ & 0.0547 & 0.0541 & 0.0537 & 0.0543 & 0.0528 & 0.0529 \\
\hline
\multirow{4}[1]{*}{1.5SUM} & $\min \varepsilon_{90}$ & 0.0262 & 0.0342 & 0.0292 & 0.0308 & 0.0314 & 0.0316 \\
      & $\max \varepsilon_{90}$ & 0.0685 & 0.0709 & 0.0713 & 0.0691 & 0.0703 & 0.0703 \\
      & ${\rm median} \varepsilon_{90}$ & 0.0617 & 0.0563 & 0.0587 & 0.0559 & 0.0556 & 0.0558 \\
      & $\varepsilon_{90}$ & 0.0553 & 0.0547 & 0.0546 & 0.0527 & 0.0531 & 0.0532 \\\hline
\multirow{4}[0]{*}{kC} & $\min \varepsilon_{90}$ & 0.0269 & 0.0368 & 0.0265 & 0.0251 & 0.0272 & 0.0272 \\
      & $\max \varepsilon_{90}$ & 0.0650 & 0.0700 & 0.0698 & 0.0709 & 0.0709 & 0.0700 \\
      & ${\rm median} \varepsilon_{90}$ & 0.0588 & 0.0564 & 0.0559 & 0.0559 & 0.0569 & 0.0571 \\
      & $\varepsilon_{90}$ & 0.0514 & 0.0549 & 0.0536 & 0.0534 & 0.0538 & 0.0535 \\\hline
\multirow{4}[1]{*}{akC} & $\min \varepsilon_{90}$ & 0.0349 & 0.0338 & 0.0360 & 0.0305 & 0.0256 & 0.0604 \\
      & $\max \varepsilon_{90}$ & 0.1042 & 0.1041 & 0.1017 & 0.3524 & 0.1100 & 0.1303 \\
      & ${\rm median} \varepsilon_{90}$ & 0.0906 & 0.0888 & 0.0820 & 0.0885 & 0.0676 & 0.0931 \\
      & $\varepsilon_{90}$ & 0.0815 & 0.0799 & 0.0778 & 0.1115 & 0.0713 & 0.0923 \\
\hline
\multirow{4}{*}{MED} & $\min \varepsilon_{90}$ & 0.0342 & 0.0329 & 0.0346 & 0.0332 & 0.0429 & 0.0270 \\
      & $\max \varepsilon_{90}$ & 0.1064 & 0.0994 & 0.0997 & 0.1102 & 0.3410 & 0.3266 \\
      & ${\rm median} \varepsilon_{90}$ & 0.0709 & 0.0872 & 0.0894 & 0.0649 & 0.0844 & 0.0714 \\
      & $\varepsilon_{90}$ & 0.0738 & 0.0784 & 0.0794 & 0.0671 & 0.1215 & 0.1012 \\
\hline
\end{tabular}}
\end{table}

\input{figura-dw3.tex}

\section{Conclusions and Further Research}
This paper introduces a new framework for fitting hyperplanes to a given set of points by considering distance-based residuals and applying generalized ordered weighted averaging aggregation criteria. Mathematical programming formulations are proposed for those models and  some properties are proven. Two important particular cases of residuals are analyzed in more detail, namely those induced by block norms or $\ell_\tau$ norms for $\tau \geq 1$. A  new goodness of fitting measure is also introduced for this framework, which extends the classical coefficient of determination in least sum of squares fitting with vertical distances. Extensive computational experiments run in Gurobi under \r{} are reported in order to illustrate and validate the new methodology for computing optimal fitting hyperplanes.

The results in this paper admit some extensions  applying similar tools. Among them we mention regularization adding constraints to overcome ill-posed data set, the simultaneous computation of several (more than one) hyperplanes to a given data set such that each single point is ``allocated''  to its \textit{closest model}. This approach would allow to analyze structural changes on the behavior of the data (in different periods of time or for different values of one of the variables). The main, non trivial, difference between those models and the ones proposed in this paper is analogous to that that exists between the so-called single-facility  and multifacility location problems (see \cite{NP05}). It is well-known that multifacility problems become easily hard even if the single-facility case were \textit{easy}. Hence, although very interesting, the above extension needs  further analysis. Another interesting extension is the use of mathematical programming tools to fit hyperplanes to binary data. The usual techniques to estimate those models are based on likelihood estimation since least squares estimation is known to get no desirable results on this type of data. Here our proposal will fit in a natural way and will deserve further attention.

%\ACKNOWLEDGMENT{The first and second authors were partially supported by the projects FQM-5849 (Junta de Andaluc\'ia $\backslash$ FEDER),  MTM2010-19576-C02-01  and MTM2013-46962-C2-1-P (MICINN, Spain).
%}

\section*{Acknowledgements}
The first and second authors were partially supported by the project MTM2016-74983-C2-1-R and MTM2013-46962-C2-1-P (MINECO, Spain).

\section*{References}

\bibliographystyle{elsarticle-harv}

%\section{References}

\end{document}